\documentclass{article}
\usepackage[utf8]{inputenc}
\usepackage{amssymb}
\usepackage{amsthm}
\usepackage{amssymb,amsmath,graphicx,epsfig,psfrag} 
\usepackage{slashed} 
\usepackage{amsmath}
\usepackage{setspace}
\usepackage{graphicx}
\usepackage{epsfig}
\usepackage{verbatim}
\usepackage{color}
\usepackage{amssymb}
\usepackage{amsthm}
\usepackage{psfrag}
\usepackage{slashed}
\usepackage{amsfonts,amsthm,bm}
\usepackage{authblk}

\def\ds{\displaystyle}
\newcommand{\bu}{{{\bf{  u}}}}
\newcommand{\bv}{{ v}}

\newcommand{\bw}{{{\bf{  w}}}}

\newtheorem{lemma}{Lemma}
\newtheorem{theorem}{Theorem}

\begin{document}
\title{Stability analysis and error estimates of 
       a projection based variational multiscale method for Oseen equations in moving domains}

\author[1]{Birupaksha Pal \thanks{Corresponding Author: birupaksha.pal@gmail.com}}
\author[2]{Sashikumaar Ganesan\thanks{sashi@iisc.ac.in}}

\affil[1,2]{Department of Computational and Data Sciences, Indian Institute of Science, Bangalore, India}

\date{May 25, 2019}
\maketitle

\begin{abstract} Stability and error estimate for the Oseen equations in a projection based variational setup has been derived in this paper.
 The use of Geometric Conservation Law (GCL) provides unconditional
stability whereas without using GCL we have a conditional scheme which imposes restriction on the time step.
Further using the stability results derived, we make the 
 first order error estimate using a backward Euler time discretization scheme.
\end{abstract}


\section*{Introduction}
  In this paper we shall discuss the mathematical analysis of the projection based variational multiscale (VMS)
      scheme that we developed for  the 
      Navier--Stokes equations (NSE) in arbitrary Lagrangian Eulerian (ALE) formulation \cite{Pal2016projection}.
      Stability and error estimates for NSE in  stationary domains have been given in a series of papers 
      by Heywood and Rannacher in the 1980s \cite{heywood1982finite,heywood1986finite,heywood1988finite,heywood1990finite}. 
      Moreover, in case of stationary domains, the stability estimate of the semidiscrete problem for a projection based VMS-NSE has been given in \cite{john2008finite2,john2008finite}.
	Fully discrete error estimates of both first and second order has been given in \cite{shang2013error}, 
	and for variants of projection based method, second order estimates has been given  in \cite{galvin2011new,rebollo2015numerical,shan2012variational}.
	Apart from these, abundant materials exist on mathematical analysis of NSE for stationary domains.
      However, when it comes to time dependent moving domains there is a real paucity of work on analysis of NSE. Stability estimate for 
      implicit Euler time discretization of ALE NSE has been presented in \cite{NOB01}, however,  to the best of the authors knowledge
      no error analysis  exists for NSE in ALE formulation. Now, ALE being one of the most popular and widely used schemes for problems with
	moving boundaries, it's analysis is imperative for a holistic study of the subject.
	
	The first steps in this direction can be considered to be the works of Nobile \cite{nobile1999stability,NOB01} 
	and Boffi, Gastaldi \cite{gastaldi2001priori,boffi2004stability},
	where the ALE formulation of non-stationary convection-diffusion scaler equations(CDE) in ALE form has been considered.
	The other important works on stability and convergence analysis
	for scaler CDE in an ALE framework can be found in \cite{formaggia2004stability}, 
	which considers an orthogonal subgrid scale stabilized(OSS) form and 
	\cite{FarhatGeuzaineGrandmont}, where the discrete geometric conservation law and nonlinear stability has been discussed.

	In this study, we consider the linear form of NSE, better known as the Oseen equations.
	  An important aspect of this study is the consideration of the projection based VMS form of the ALE-NSE,
	  where the effect of the subgrid scale model is accounted for as added diffusion.
	  The subgrid scale model term in the VMS form we developed acts as  additional viscosity and is agreeable to regular treatment as that of the diffusion term in the NSE.
	If we approach the problem for the Oseen equations in a manner similar to the scaler CDE, then the mesh velocity 
	  in the convective part gives rise to some extra terms, which gets quite challenging to bound.
	In this work our strategy will be to first derive a stability estimate of the semidiscrete 
	form of  ALE-Oseen equations in a projection based  VMS frame, in the footsteps of \cite{NOB01,boffi2004stability},
	and finally this stability estimate shall be used to derive an
	error estimate for the fully discrete system. We consider implicit Euler time discretization for the fully discrete form.
	
      
      In the next  sections  we shall describe the Oseen equations, define the relevant spaces to derive it's variational form, 
      and  get to it's time discretized form. In the third section, we discuss some preliminary but important results
      which shall be put to later use. 
      Then in the next two subsequent sections we derive a
      stability analysis of the Oseen equations. We derive two estimates one using the geometric conservation law which gives us an
      estimate independent of the domain deformation and the other a more generalized result without taking recourse
      to GCL, this leads us to an estimate which is dependent on the domain deformation or the mesh velocity.
      However, by a smart choice of time step length the estimate can be considered to be independent of 
      mesh movement for all practical purpose. Finally, we derive an estimate for error due to the time discretization
     by implicit Euler time scheme, in obtaining the final estimate we make use of the more generalized
      stability estimate without GCL.
\section{Navier-Stokes equations for time dependent domains}
We consider an incompressible fluid flow in a moving/deforming domain which is described by the time-dependent incompressible Navier--Stokes equations:
\begin{equation}\label{model}
\begin{array}{rcll}
\ds\frac{\partial\bu}{\partial t}-2\mu \nabla \cdot \mathbb{D}(\bu)
+ (\bu\cdot \nabla) \bu +\nabla p &=& \bf{f}  & \text{ in} \,\ (0,\rm{T}] \times \Omega_t,\\
\nabla \cdot \bu &=& 0 & \text{ in} \,\ (0,\rm{T}] \times \Omega_t. 
\end{array}
\end{equation}
 Here, $\bu = (u_1,u_2,u_3)^{{T}}$ is the fluid velocity,  $p$ is the pressure in the fluid, $\text{T}$ is a given final time,
   and $\Omega_t\subset \mathbb{R}^{3}$, $t\in(0,\text{T}]$ is a time-dependent domian.
   The NSE are closed with the initial condition
\begin{align*}
 \bu(0,\cdot) &= \bu_{0} ~\text{ in }~\Omega_0 \\
\end{align*}

and boundary conditions
$$
 \textbf{u} = \textbf{g}\text{(t)}  ~\text{on}~  \partial\Omega_t 
$$
 Here, $\textbf{g}\text{(t)} $ is a continous function denoting the domain deformation, $\bu_0$ is a given initial velocity and $\bu_D= (\bu_{in},0)^T$ a given inlet velocity, 
$\mathbb{I}$ is the identity tensor and $\textbf{n}$ is the outward normal to the boundary $\Gamma_{out}$.
Further, the velocity deformation tensor is defined as
$
\mathbb{D}(\bu) = \displaystyle\frac{\nabla \bu +\nabla \bu^{T} }{2}, \quad \text{and} ~\mu ~ \text{represents the inverse of Reynolds number}.
$


%
%
%
%
%
%
%

\subsection{ALE formulation}
\label{sec:1}

In order to handle the time-dependent domain, we now derive an arbitrary Lagrangian-Eulerian form of  NSE~\eqref{model}. 
Let $\hat\Omega \subset \mathbb{R}^3$ be a reference domain, and define a family of bijective ALE mappings
\[
\mathcal{A}_t:\hat{\Omega} \rightarrow \Omega_t,  \qquad  \mathcal{A}_t(Y)= x(t,Y), \qquad t\in(0,\rm{T}).
\]
The reference domain  $\hat\Omega$ can simply be the initial domain $\Omega_0$ or the previous time-step domain when the deformation in the domain is large.
In addition, for a  scaler function $v\in C^0({{\Omega_t}})$ on the Eulerian frame, we define it's corresponding function $\hat v\in C^0({{\hat\Omega} })$ on the ALE frame as 
\[
\hat v : (0, {\rm{T}}) \times \hat\Omega  \rightarrow \mathbb{R}, \qquad 
 \hat{v}:=v\circ \mathcal{A}_t,\quad \text{with} \quad \hat{v}(t,Y) = v(t,\mathcal{A}_t(Y)).
\]
Further, the time derivative   on the ALE frame is defined by 
\[
\ds\frac{\partial v}{\partial t} \Big|_{Y}:(0, {\rm{T}}) \times \Omega_t  \rightarrow \mathbb{R},~ 
\ds\frac{\partial v  }{\partial t} \Big|_{Y }(t,x) = \ds\frac{\partial \hat v }{\partial t}(t,Y),\quad \text{with} \quad Y= \mathcal{A}_t^{-1}(x).
 \]
We now apply the chain rule to the time derivative of $v\circ\mathcal{A}_t$ on the ALE frame to get
\begin{align*}
\frac{\partial v}{\partial t} \Big|_{Y} &= \ds\frac{\partial v}{\partial t} (t,x) + \ds\frac{\partial x}{\partial t}\Big|_{Y}\cdot\nabla_x v  = 
 \ds\frac{\partial v}{\partial t} + \ds\frac{\partial \mathcal{A}_t(Y)}{\partial t}  \cdot\nabla_x v = \ds\frac{\partial v}{\partial t} + \bw\cdot\nabla_x v, 
\end{align*}
where $\nabla_x$ represents the divergence function and  $\bw$ is the domain velocity. 
Using similar arguements to NSE~\eqref{model} to account for the  deformation in the domain,   we get the ALE form of the NSE as
\begin{equation}\label{ALEmodel}
\begin{array}{rcll}
\ds\frac{\partial\textbf{u}}{\partial t}\Big|_{Y}-\frac{2}{\text{Re}}\nabla \cdot \mathbb{D}(\textbf{u}) + ((\textbf{u}-\textbf{w})\cdot \nabla \textbf{u}) + \nabla p &=& \bf{f};~~~
\nabla \cdot \bu = 0. 
\end{array}
\end{equation}
Note that the main difference between equations ~\eqref{model} and~\eqref{ALEmodel} is the additional domain velocity $\bw$ in the ALE form   that accounts for the deformation of the domain.
The ALE form~\eqref{ALEmodel} can be viewed as a generalized form of NSE, since the Lagrangian form of NSE can be obtained by
setting $\bw = \bu$ and the Eulerian form of NSE can be obtained by setting $\bw = 0$.
\subsection{Linearized form of the ALE-NSE or ALE Oseen equations}
\label{sec:2}

The ALE-NSE equations derived in the previous section have an associated complex non-linearity in the convective acceleration term.
However, for all the practical purposes of simulation and numerical computation of fluid flow in time-dependent domians 
this convection term is linearized using methods like fixed point iteration among other such linearization techniques. 
Linearization of the convective term in  ALE-NSE leads to the Oseen equations which are as follows: 

\begin{align}\label{OseenEqns}
\ds\frac{\partial\textbf{u}}{\partial t}-2\text{Re}^{-1}\nabla \cdot \mathbb{D}(\textbf{u}) + ((\textbf{u}^*-\textbf{w})\cdot \nabla) \textbf{u} 
+ \nabla p &= \textbf{f} & \text{ in} \,\ (0,\text{T}] \times \Omega_t\\ \nonumber
\nabla \cdot \textbf{u} &= 0   & \text{ in} \,\ (0,\text{T}] \times \Omega_t
\end{align}

\text{with the boundary conditions}

\begin{align*}
\textbf{u}(0,\cdot) &= \textbf{u}_{0}  ~~~\text{ in} \,\ \Omega_0\\ 
\textbf{u} &= \textbf{g}\text{(t)}  ~~~\text{ on} \,\ \partial\Omega_t 
\end{align*}
Moreover, the source term $\textbf{ f} \in [\rm{H}^{-1}(\Omega_{t})]^3$ and we consider $\bu^* \in [W^{1,\infty}(\Omega_t)]^3$ 
such that $(\nabla \cdotp \bu^*,q)=0 ~\text{for all}~ q \in \rm{L}^2(\Omega_t)$.

Let us define the following spaces :
 $$
      \begin{array}{rcll}
     
     
      \hat{Q}(\Omega_0) &=& \{ q : q \in \rm{L}^2(\Omega_0) , \int_{\Omega_0} q dx =0 \} \\
      
       \hat{X}(\Omega_0) &=& \{ \textbf{v} : \textbf{v} \in (\rm{H}^1(\Omega_0))^3 ~~\text{and}~~ (\nabla \cdotp \textbf{v},q) = 0  ~~\forall q \in \hat{Q}(\Omega_0) \} \\

      \hat{V}(\Omega_0) &=& \{ \textbf{v} \in \hat{X}(\Omega_0) :  \textbf{v}|_{\partial\Omega_0} = 0   \} \\ 

      X(\Omega_t) &=& \{ \textbf{v} : (0,\text{T}]\times \Omega_t \rightarrow \mathbb{R}^3 ,\textbf{v}= \hat{\textbf{v}} \circ {A_t}^{-1} , \hat{\textbf{v}} \in \hat{X}({\Omega_0}) \} \\

      Q(\Omega_t) &=& \{ q : (0,\text{T}]\times \Omega_t \rightarrow \mathbb{R} ,q= \hat{q} \circ {A_t}^{-1} , \hat{q} \in \hat{Q}({\Omega_0}) \} \\
      
      V(\Omega_t) &=& \{ \textbf{v} : (0,\text{T}]\times \Omega_t \rightarrow \mathbb{R}^3 ,\textbf{v}= \hat{\textbf{v}} \circ {A_t}^{-1} , \hat{\textbf{v}} \in \hat{V}(\Omega_0) \} \\
      \end{array}
      $$

\noindent As we are dealing with functions with both  time and space variables  we shall be using the following notations:
$$
\rm{L}^2(\text{I}; \rm{H}^p(\Omega)) =\{ \textbf{v}: \text{I}\rightarrow \rm{H}^p(\Omega) | \int_\text{I} \parallel \textbf{v}(t) \parallel^2_{\rm{H}^p(\Omega)} dt < \infty \},
$$
which is equipped with the norm 
$$
\parallel \textbf{v} \parallel_{\rm{L}^2(\text{I}; \rm{H}^p(\Omega))} = \left( \int_\text{I} \parallel \textbf{v}(t) \parallel^2_{\rm{H}^p(\Omega)} dt \right)^\frac{1}{2}
$$
Here, $\text{I}$ denotes the entire time interval $(0, \text{T}]$.

\noindent Multiplying \eqref{OseenEqns} with velocity test function from  $V(\Omega_t)$ and pressure test function from $Q(\Omega_t)$,
      we get:\\ \noindent Find $\textbf{u} ~\in~ X(\Omega_t)~$ such that the 
      following holds for all $\textbf{v}~~\in~ V(\Omega_t)$ 

      \begin{equation}\label{aleoseen}
      \ds \left( \frac{\partial \textbf{u}}{\partial t} ,\textbf{v} \right) + 2\mu (\nabla \textbf{u}, \nabla \textbf{v}) + (((\textbf{u}^* -\textbf{w})  \cdotp \nabla) \textbf{u} ,\textbf{v}) =(\bf f,\textbf{v}) 
      \end{equation}

\subsection{Linearized ALE-NSE-VMS}
\label{sec:3}
\noindent We can use the following notation for denoting \eqref{aleoseen}:

\begin{align}\label{weakALE}
A(\textbf{u} ;\textbf{u};\textbf{v}) - b( \bw, \textbf{u},\textbf{v})=0 
\end{align} 
where
 \begin{align*}
 A(\bu;\textbf{v}))    &= \left(\ds\frac{\partial\textbf{u}}{\partial t}\Big|_{Y} , \textbf{v}\right) + 
 b(\bu^*,\textbf{u},\textbf{v}) \nonumber   +\displaystyle{2\mu} \left(\nabla(\textbf{u}),\,\nabla(\textbf{v})\right)\\
  b(\bw,\textbf{u},\textbf{v}) &= (\bw \cdot \nabla\textbf{u},\textbf{v}).
 \end{align*}
Here,  $(\cdot,~\cdot)$ denotes the $\rm{L}^2-$inner product in   $\Omega_t$.
Further, $\bu\in \rm{L}^2(0,\text{T};X(\Omega_t))$   implies that the mapping $t\mapsto \bu(t)$ is continuous.  

The solution of the variational form \eqref{weakALE} can be evaluated numerically using a standard finite element approach. 
However if the Reynolds number of the flow is high, the equation \eqref{weakALE} is representative of turbulent flow.
Now, as flows with high Reynolds number tend to exhibit multitude of scales in it's flow fields, 
use of standard Galerkin method would not be able to capture the fluid flow scales 
that are smaller than the mesh size (discretization parameter).
Moreover, if  the mesh were to be taken fine enough to capture the dissipation of energy at the
Kolmogorov length scale it would require enormous memory and computing power.
This lends the problem to be a perfect candidate for multiscale methods.
It's basic idea is to decompose the flow fields into resolved (large) and unresolved (small)
scales and incorporate the effects of small-scales into the solution of the large-scales by a turbulence model.
We choose a reasonable mesh size to capture the large scale flow dynamics, whereas the unresolved or small scale flow dynamics 
is approximated by a turbulence model. This approximation of the small unresolved scales thus corrects the discretization error
caused due to the limitation of the mesh chosen to represent the large scales.
 A simple and effective approach for it is projection based variational multiscale method. 
A standard finite element space can be used to represent the resolved scales. 
Whereas, the remnant of the solution $i.e,$ the small scales also known as the sub-grid scales
is infinite dimensional and is modeled.

 \noindent In this paper,we shall use a three-scale decomposition, where 
the resolved solution space is  decomposed  into resolved large and resolved small scales, $i.e,$
\begin{equation*}\label{soln_space_decomp}
 V = \overline{V} \oplus \tilde{V} \oplus \widehat{V} \quad \text{and the velocity test space as} \quad V_0 = \overline{V_0} \oplus \tilde{V_0} \oplus \widehat{V_0}.
\end{equation*}

\noindent  the above equations, the  bar, the tilde and the hat over the spaces represent the resolved large, the resolved small and the  unresolved small scales, respectively.
Consequentially, the functions $\bu\in V$ can be written as
\begin{equation}\label{scale_decomp}
 \textbf{u} = \overline{\textbf{u}} + \tilde{\textbf{u}} + \widehat{\textbf{u}}.
\end{equation}

Using the decomposition \eqref{scale_decomp},  we can write the momentum balance equation \eqref{weakALE} as:  

\begin{align}
 A(\textbf{u}^*;\overline{\textbf{u}},\overline{\textbf{v}})  + A(\textbf{u}^*; \tilde{\textbf{u}},\overline{\textbf{v}}) +
 A(\textbf{u}^*;\widehat{\textbf{u}},\overline{\textbf{v}}) -
 b(\bw,\overline{\textbf{u}} + \tilde{\textbf{u}} + \widehat{\textbf{u}}, \overline{\textbf{v}} )&=0 \nonumber \\
  A(\textbf{u}^*;\overline{\textbf{u}},\tilde{\textbf{v}})+ A(\textbf{u}^*;\tilde{\textbf{u}},\tilde{\textbf{v}})  +
  A(\textbf{u}^*;\widehat{\textbf{u}},\tilde{\textbf{v}})-
  b(\bw,\overline{\textbf{u}} + \tilde{\textbf{u}} + \widehat{\textbf{u}}, \tilde{\textbf{v}} ) &=0 \label{3scaleVMS}  \\
  A(\textbf{u}^*;\overline{\textbf{u}},\widehat{\textbf{v}})+A(\textbf{u}^*;\tilde{\textbf{u}},\widehat{\textbf{v}})  +
  A(\textbf{u}^*;\widehat{\textbf{u}},\widehat{\textbf{v}})-
  b(\bw,\overline{\textbf{u}} + \tilde{\textbf{u}} + \widehat{\textbf{u}}, \widehat{\textbf{v}} )&=0. \nonumber
\end{align}

\emph{Modeling assumptions}\\
\indent The following assumptions are made on \eqref{3scaleVMS}
\begin{itemize}
  \item The last equation \eqref{3scaleVMS}  consists of test functions from the space of the unresolved scales and is thus ignored
  \item The first equation in \eqref{3scaleVMS} is the one with the test function from the resolved large scales. 
 In this equation the last term $b(\bw,\overline{\textbf{u}} +  \tilde{\textbf{u}} + \widehat{\textbf{u}}, \overline{\textbf{v}} )$
 which contains the mesh velocity $\bw$ incorporates the information about the domain movement into the model. 
 We can expand this term linearly and then among the terms we get by expansion $b(\bw,\overline{\textbf{u}}, \overline{\textbf{v}} )$,
 $b(\bw, \tilde{\textbf{u}}, \overline{\textbf{v}} )$ can be combined with
 $A(\textbf{u}^*;\overline{\textbf{u}},\overline{\textbf{v}})$
 and $A(\textbf{u}^*;\textbf{u}; \tilde{\textbf{u}},\overline{\textbf{v}})$, respectively, to give
 $A(\textbf{u}^*-\bw;\overline{\textbf{u}},\overline{\textbf{v}})$ and $A(\textbf{u}^*-\bw; \tilde{\textbf{u}},\overline{\textbf{v}})$.
 We can also modify the second equation in \eqref{3scaleVMS} similarly.
 \item  The Reynolds stress and cross stress terms  
 \begin{align*}
  A( {\textbf{u}^*}-\textbf{w} ;\widehat{\textbf{u}},\overline{\textbf{v}}), \,   
     b(\widehat{\textbf{u}},  \overline{\textbf{u}},  \overline{\textbf{v}} ), \, b(\widehat{\textbf{u}},  
     \tilde{\textbf{u}},  \overline{\textbf{v}} ),  \, b(\widehat{\textbf{u}},  \widehat{\textbf{u}},  \overline{\textbf{v}} ) 
 \end{align*}
that contain the unresolved and resolved large scales  are assumed to be zero and it is further  assumed that the direct influence
 of the unresolved small scales  on the resolved small scale is zero.
 \item The influence of unresolved scales on resolved small scales is modeled by an turbulence model, that is,
\begin{align*}
 A(\textbf{u}^*-\textbf{w};\widehat{\textbf{u}},\tilde{\textbf{v}}) + b(\widehat{\textbf{u}} ,  \overline{\textbf{u}},  
 \tilde{\textbf{v}} ) + b(\widehat{\textbf{u}} ,  \tilde{\textbf{u}},  \tilde{\textbf{v}})  + b(\widehat{\textbf{u}},  \widehat {\textbf{u}}, 
 \tilde{\textbf{v}} ) \approx B(\textbf{u}^*-\textbf{w};\overline{\textbf{u}},\tilde{\textbf{u}},\tilde{\textbf{v}}).
\end{align*}
\end{itemize}
 This modeling is essential to incorporate the effect of the unresolved scales into the resolved scales of the flow. 
Imposing these assumptions, VMS form of   ALE-NSE \eqref{3scaleVMS} reads
\begin{align}\label{3_VMS}
 A(\bu^*-\textbf{w};\overline{\textbf{u}},\overline{\textbf{v}})   
 + A(\bu^*-\textbf{w};\tilde{\textbf{u}},\overline{\textbf{v}}) & = 0\nonumber\\
 A(\bu^*-\textbf{w};\overline{\textbf{u}},\tilde{\textbf{v}})     
 + A(\bu^*-\textbf{w};\tilde{\textbf{u}},\tilde{\textbf{v}}) 
 + B(\bu^*-\textbf{w};\overline{\textbf{u}},\tilde{\textbf{u}},\tilde{\textbf{v}}) &= 0 \nonumber
\end{align}
Due to these modifications it can  clearly be seen that the turbulent model in the three-scale 
VMS acts only on the resolved small scales and not on the resolved large scales directly.
Nevertheless the three-scale model  incorporates the effects of unresolved scales in the resolved large scales  indirectly through
the coupling of  resolved small scales with  resolved large scales. Though VMS is a variant of LES, VMS differs fundamentally 
 from the traditional LES due to the fact that the 
turbulent model in LES acts directly on all resolved scales.

\vspace{5mm}
\subsection{Discretization of the Oseen equations in ALE-VMS form }
      In a finite element based approach using triangulation we discretize the domain $\Omega_t$ into  $\Omega_{t,h}$.
      On this discretized domain we consider the inf-sup stable
      finite dimensional subspaces $X_h(\Omega_{t,h}), Q_h(\Omega_{t,h})$ of $X(\Omega_t),Q(\Omega_t)$ respectively.
      We further define the space of discretely divergence free functions $$ V_h(\Omega_{t,h}) = \{ \textbf{v}_h \in X_h(\Omega_{t,h}) | (\nabla \cdotp \textbf{v}^h,q^h) = 0  ~~\forall q^h \in Q_h(\Omega_{t,h})  \},$$
      but for the sake of convenience shall denote it as $V_h(\Omega_{t})$.
      Hence, the semi discrete form of the Oseen problem in ALE-VMS form becomes :\\
      Find $\textbf{u}_h= \overline{\textbf{u}}_h + \tilde{\textbf{u}}_h  ~\in~ V_h(\Omega_t) $ such that the following holds 
      for all $\textbf{v}_h= \overline{\textbf{v}}_h + \tilde{\textbf{v}}_h ~\in~ V_h(\Omega_t) \cap V(\Omega_t), $ 

      \begin{equation}\label{aleoseen_semidiscrete}
       \ds \left(  \frac{\partial \textbf{u}_h}{\partial t} ,\textbf{v}_h \right) + 2\mu (\nabla \textbf{u}_h, \nabla \textbf{v}_h) + \mu_T (\nabla \tilde{\textbf{u}}_h, \nabla \tilde{\textbf{v}}_h)+
      (((\textbf{u}^* -\textbf{w}_h)  \cdotp \nabla) \textbf{u}_h ,\textbf{v}_h) =({\bf f},\textbf{v}_h)
      \end{equation}
      Here,  $\mu_T$ represents the additional turbulent viscosity
       added to the resolved small scales which is further considered to be constant.
      

\subsection{Discretization of the ALE mapping}
Let the considered time interval $[0,  {\rm T}]$ be decomposed of  into $N$ equal subintervals $0=t^0<t^1<\dots <t^N={T}$, 
 and we further denote the uniform time step by $\Delta t   $ = $t^{n}$ - $t^{n-1}$, $1\le n \le N$. 
Let   $\bu_h^n$ be an
approximation of $\bu(t^n,x)$ in $V_{h}(\Omega_{t^n})$, where $\Omega_{t^n}$ is the deforming domain at time $t=t^n$. 
We next define the semi-discrete mesh velocity $\bw_h$ in space using the semi-discrete ALE mapping
\begin{equation}\label{discALE}
 \mathcal{A}_{h,t}:\hat{\Omega}_h \rightarrow \Omega_{h,t}.
\end{equation}
We first discretize  the ALE mapping in time using a linear interpolation. Denoting the discrete ALE mapping by 
$\mathcal{A}_{h, \Delta t}$, we define it for every $\tau\in[t^n,t^{n+1}]$ by

\begin{equation}\label{Chap4DiscALE}
 \mathcal{A}_{h, \Delta t}(Y) = \frac{\tau-t^n}{\Delta t}\mathcal{A}_{h,t^{n+1}}(Y) +  \frac{t^{n+1}-\tau }{\Delta t}\mathcal{A}_{h,t^{n}}(Y),
\end{equation} 

where $\mathcal{A}_{h,t}(Y)$ is the time continuous ALE mapping defined in~\eqref{discALE}. Since the discrete ALE mapping is defined linearly in time, we obtain the discrete mesh velocity 
\begin{equation}\label{Chap4Discmeshvelo}
 \hat\bw_h^{n+1}(Y)= \frac{\mathcal{A}_{h,t^{n+1}}(Y) - \mathcal{A}_{h,t^{n}}(Y)}{\Delta t}
\end{equation}
as a piecewise constant function in time. Further, we define the mesh velocity on the Eulerian frame as
\[
  \bw_h^{n+1} =  \hat\bw_h^{n+1}  \circ\mathcal{A}^{-1}_{h,\Delta t}(x).
\]

The discrete ALE mapping lies in the space $[W^{1,\infty}(\Omega_0 \times I)]^d$
, and satisfies the hypothesis of proposition 1.3.1 of  ~\cite{NOB01}. We shall later use the property of the mapping namely 
$\parallel \mathcal{A}_{h,t} \parallel_{L^{\infty}(\Omega_0 \times I)}$ to be bounded in our error estimates.
 \section{Preliminary results}
      
      Before moving onto the stability and error estimates, we recount some important concepts and results
      that will be used in our analysis.
	      
      \begin{lemma}(Gronwall lemma,  \cite{heywood1982finite})
	Let   $\Delta t,~\gamma_n,~f,~A_n,~ B_n,~ C_n $ be given sequences of non-negative numbers for $n \geq 0$ such that the following inequality holds
	\[
	A_n + \Delta t \sum_{i=0}^{n} B_i \leq \Delta t  \sum_{i=0}^{n} \gamma_i A_i + \Delta t \sum_{i=0}^{n} C_i + f. 
	\]
	Suppose that $\Delta t \gamma_i   < 1$ for all $i$, and set $\sigma_i = \left(1-\Delta t \gamma_i  \right)^{-1}$. We then have
	\[
	A_n + \Delta t \sum_{i=0}^{n} B_i \leq \exp \left(\Delta t \sum_{i=0}^{n} \sigma_i \gamma_i\right)\left[\Delta t \sum_{i=0}^{n} C_i + f\right].
	\]
	\end{lemma}
	\subsection{Reynolds transport theorem}\label{ReyTransTh}
	It is one of the most useful results for studying PDEs in time-dependent domains, especially in ALE frame. It states that
	for any $\varphi( t, {\bf x} ): \Omega_t \rightarrow \mathbb{R}$ such that $\varphi = \hat \varphi~ \circ~ \mathcal{A}_t^{-1} $ for some
	$\varphi: \Omega_0 \rightarrow \mathbb{R}$. the following equation holds:	
	\begin{equation}
	 \frac{d}{dt} \int_{\Omega_{t}} \varphi ~\text{d}\Omega = \int_{\Omega_{t}} \varphi (\nabla \cdotp \bw)  ~\text{d}\Omega.
	\end{equation}

      \subsection{Geometric conservation law (GCL)}\label{GCLcond0}
      The Geometric conservation law (GCL) was first introduced in \cite{thomas1979geometric} and later studied closely in~\cite{farhat2001discrete,guillard2000significance} as a minimum criterion for stability.
      It is a sufficient condition which ensures that a time stepping scheme will be accurate  in time upto the first order for a moving domain.
      However, it does not guarantee the order of accuracy of the time-integration. It has been observed 
      that apart from the first order time stepping schemes, GCL is neither a necessary nor a sufficient condition for achieving the expected temporal accuracy \cite{NOB01}.
      
      The GCL can be formulated as time integration of the term containing the mesh velocity in the ALE form of NSE.
      In the finite element context it can be written as 
      \begin{equation}\label{GCLcond1}
       \int_{\Omega_{t^{n+1}}} \bm{\phi_i} \cdotp \bm{\phi_j} ~\text{d}\Omega_{t^{n+1}} - \int_{\Omega_{t^{n}}} \bm{\phi_i} \cdotp \bm{\phi_j} ~\text{d}\Omega_{t^{n}}
       = \int_{t^n}^{t^{n+1}}\int_{\Omega_{t}} \bm{\phi_i} \cdotp \bm{\phi_j} (\nabla \cdotp \bw_h) \, \text{d}\Omega_t \,dt.
      \end{equation}
      Here, $\bm\phi_i, \bm\phi_j\in V_h(\Omega_{t}). $
      This relation can be derived using Reynolds transport theorem. If $\bm\phi_i$ is taken equal to $\bm\phi_j$,
      then it reduces to the following form. 
      \begin{equation*}\label{GCLcond2}
       \int_{\Omega_{t^{n+1}}} \vert \bm{\phi_i} \vert^2 ~\text{d}\Omega_{t^{n+1}} - \int_{\Omega_{t^{n}}} \vert \bm{\phi_i} \vert^2 ~\text{d}\Omega_{t^{n}}
       = \int_{t^n}^{t^{n+1}}\int_{\Omega_{t}} \vert \bm{\phi_i} \vert^2 (\nabla \cdotp \bw_h) ~\text{d}\Omega_t dt.
      \end{equation*}
      
      \section{Stability analysis for modified implicit Euler time discretization applied to conservative form of ALE-VMS Oseen equations}
      The ALE formulation involves mapping between different domain configurations. For the ease of 
      notations further down wherever $\textbf{u}^n_h$ is defined on domain $\Omega_p$ 
      it implicitly means that a transformation mapping  $\mathcal{A}_{t^n,p}= \mathcal{A}_{h,p} \circ \mathcal{A}^{-1}_{h,t^n}$ has been applied to it ,$i.e.$
      \[
       \int_{\Omega_p} \textbf{u}^n_h ~\text{d}\Omega_p =  \int_{\Omega_p} \mathcal{A}_{t^n,p} \circ \textbf{u}^n_h ~\text{d}\Omega_p.
      \]
      We shall use the standard $L^2$ norm defined as follows:
      \[
      \parallel \bu^n_h \parallel_{n} = \left( \int_{\Omega_{t^n}} \bu^n_h \cdotp \bu^n_h ~\text{d}\Omega_{t^n}\right)^{\frac{1}{2}} = \left( \int_{\Omega_{t^n}} \vert \bu^n_h \vert^2 ~\text{d}\Omega_{t^n}\right)^{\frac{1}{2}}.
      \]
       \noindent  We consider the discrete ALE mapping $\mathcal{A}_{h,t} \in [ W^{1,\infty} (I \times \Omega_0)]^d $. Further,
       we denote by $J_{\mathcal{A}^t}$ the matrix of $\mathcal{A}_{h,t}$ and for further simplicity in notations we will use $~\text{d}\Omega$ for all integrations.\\
      The conservative formulation of the semi-discrete form \eqref{aleoseen_semidiscrete} reads as follows:
      \begin{equation} \label{semidiscreteALEVMS}
      \begin{array}{rcll}
      \ds {\frac{d}{dt}}  \int_{\Omega_t} \textbf{u}_h \cdotp \textbf{v}_h ~\text{d}\Omega + 2\mu\int_{\Omega_t}  \nabla \textbf{u}_h  : \nabla \textbf{v}_h~\text{d}\Omega + 
      \mu_T\int_{\Omega_t} \nabla \tilde{\textbf{u}}_h : \nabla \tilde{\textbf{v}}_h~\text{d}\Omega \\
      + \int_{\Omega_t} \nabla \cdotp [(\textbf{u}^*- \textbf{w}_h )\otimes \textbf{u}_h ] \cdotp \textbf{v}_h~\text{d}\Omega = \int_{\Omega_t} {\bf f} \cdotp \textbf{v}_h~\text{d}\Omega
      \end{array}
      \end{equation}
      
      Next, we consider the  implicit Euler time discretization scheme to the semi discrete form \eqref{semidiscreteALEVMS}. In this scheme we consider a midpoint rule for the quadrature formula,  
      as we have considered a linear deformation of the domain. 
      Moreover, as implicit Euler is a first order accurate scheme, hence a midpoint rule for 
       time discretization satisfies the GCL. It should be noted here that instead of this midpoint rule one 
      can also choose the classical definition of backward  Euler scheme, but GCL would not be satisfied then. This time discretization
      leads us to the following estimate.
      
      \begin{lemma} \label{lemma1}
      Let ${\textbf{u}^{n}_h} \in V_{h}(\Omega_{t^n})$ be as defined above and the regularity assumptions
      $\textbf{ f} \in [H^{-1}(\Omega_{t})]^d$ and $\textbf{u}^{0}_h \in (L^2(\Omega))^d$ be satisfied,
      then the following stability estimate holds for all $t \in [0,T]$
    \begin{align*}   
    \ds &\int_{\Omega_{t^{n+1}}} |{\textbf{u}^{n+1}_h}|^2  ~\text{d}\Omega
    +   \triangle t \sum\limits_{i=0}^n \Bigg( 3\mu \int_{\Omega_{t^{i+\frac{1}{2}}}} |\nabla {\textbf{u}^{i+1}_h}|^2 ~\text{d}\Omega+ 
    \mu_T\int_{\Omega_{t^{i+\frac{1}{2}}}}  |\nabla {\tilde{\textbf{u}}^{i+1}_h}|^2 ~\text{d}\Omega\Bigg)\\ \nonumber
    &\leqslant  {{\parallel {\textbf{u}^{0}_h} \parallel}^2}_0 +
    \triangle t \Big( \frac{1+C_\Omega}{\mu} \Big)\sum\limits_{i=0}^n  {{\parallel {\bf f}^{i+\frac{1}{2}} \parallel}^2}_{H^{-1}(\Omega_{t^{i+\frac{1}{2}}})}    
    \end{align*}
    \begin{proof}

      The fully discrete form of \eqref{semidiscreteALEVMS} using modified Euler time discretization scheme reads:
      \begin{align}
      \ds  &\int_{\Omega_{t^{n+1}}} {\textbf{u}^{n+1}_h} \cdotp \textbf{v}_h~\text{d}\Omega  -  \int_{\Omega_{t^{n}}} {\textbf{u}^{n}_h} \cdotp \textbf{v}_h~\text{d}\Omega 
      +  2\mu \triangle t  \int_{\Omega_{t^{n+\frac{1}{2}}}} \nabla {\textbf{u}^{n+1}_h} : \nabla \textbf{v}_h~\text{d}\Omega  \nonumber \\ \nonumber
      &+ \mu_T \triangle t  \int_{\Omega_{t^{n+\frac{1}{2}}}} \nabla \tilde{{\textbf{u}}}^{n+1}_h : \nabla \tilde{\textbf{v}}_h~\text{d}\Omega + \label{ImpEulerALENSEVMS}
      \triangle t \int_{\Omega_{t^{n+\frac{1}{2}}}} \nabla \cdotp [(\textbf{u}^* - \textbf{w}_h )\otimes {\textbf{u}^{n+\frac{1}{2}}_h} ] \cdotp \textbf{v}_h~\text{d}\Omega \\&= 
      \triangle t \int_{\Omega_{t^{n+\frac{1}{2}}}} {\bf f}^{n+\frac{1}{2}}\cdotp \textbf{v}_h~\text{d}\Omega
      \end{align}

      Now in the previous equation take $\textbf{v}_h = {\textbf{u}^{n+1}_h} $, this leads to the following:

      \begin{align}
      &\ds  \int_{\Omega_{t^{n+1}}} |{\textbf{u}^{n+1}_h}|^2~\text{d}\Omega - \ds  \int_{\Omega_{t^{n}}} {\textbf{u}^{n}_h} \cdotp {\textbf{u}_h}^{n+1} ~\text{d}\Omega
      + 2\mu \triangle t \int_{\Omega_{t^{n+\frac{1}{2}}}} |\nabla {\textbf{u}^{n+1}_h}|^2 ~\text{d}\Omega \nonumber \\
      &+\mu_T \triangle t \int_{\Omega_{t^{n+\frac{1}{2}}}} |\nabla \tilde{{\textbf{u}}^{n+1}_h}|^2 ~\text{d}\Omega+
      \triangle t \int_{\Omega_{t^{n+\frac{1}{2}}}} (\nabla \cdotp (\textbf{u}^* - \textbf{w}_h )){\textbf{u}^{n+1}_h} \cdotp{\textbf{u}^{n+1}_h}~\text{d}\Omega \nonumber \\ 
      &+ \triangle t \int_{\Omega_{t^{n+\frac{1}{2}}}} [((\textbf{u}^* - \textbf{w}_h ) \cdotp \nabla ){\textbf{u}^{n+1}_h}] \cdotp {\textbf{u}^{n+1}_h}~\text{d}\Omega\\\nonumber
      &= \triangle t  \int_{\Omega_{t^{n+\frac{1}{2}}}} {\bf f}^{n+\frac{1}{2}} \cdotp {\textbf{u}^{n+1}_h}~\text{d}\Omega
      \end{align}      
      Which on simplification reads,
      \begin{align}
       &\ds  \int_{\Omega_{t^{n+1}}} |{\textbf{u}^{n+1}_h}|^2~\text{d}\Omega -  \int_{\Omega_{t^{n}}} {\textbf{u}^{n}_h} \cdotp {\textbf{u}^{n+1}_h}~\text{d}\Omega 
      + 2\mu \triangle t \int_{\Omega_{t^{n+\frac{1}{2}}}} |\nabla {\textbf{u}^{n+1}_h}|^2~\text{d}\Omega \nonumber\\\nonumber
      &+ \mu_T \triangle t \int_{\Omega_{t^{n+\frac{1}{2}}}} |\nabla \tilde{{\textbf{u}}^{n+1}_h}|^2~\text{d}\Omega + 
      \triangle t \int_{\Omega_{t^{n+\frac{1}{2}}}}( \nabla \cdotp (\textbf{u}^* - \textbf{w}_h ))|{\textbf{u}^{n+1}_h}|^2~\text{d}\Omega  \\ 
      &+ \frac{\triangle t}{2}  \int_{\Omega_{t^{n+\frac{1}{2}}}} (\textbf{u}^* - \textbf{w}_h ) \cdotp \nabla |{\textbf{u}^{n+1}_h}|^2~\text{d}\Omega  \\ \nonumber
      &= \triangle t \int_{\Omega_{t^{n+\frac{1}{2}}}} {\bf f}^{n+\frac{1}{2}}\cdotp{\textbf{u}^{n+1}_h}~\text{d}\Omega 
      \end{align}
      Now using integration by parts and combining the integrals associated with the mesh velocity we get
      \begin{align}
       &\ds  \int_{\Omega_{t^{n+1}}} |{\textbf{u}^{n+1}_h}|^2 ~\text{d}\Omega-   \int_{\Omega_{t^{n}}} {\textbf{u}^{n}_h} \cdotp {\textbf{u}^{n+1}_h} ~\text{d}\Omega
      + 2\mu \triangle t \int_{\Omega_{t^{n+\frac{1}{2}}}} |\nabla {\textbf{u}^{n+1}_h}|^2~\text{d}\Omega \nonumber \\ \nonumber
      &+ \mu_T \triangle t \int_{\Omega_{t^{n+\frac{1}{2}}}} |\nabla {\tilde{\textbf{u}}^{n+1}_h}|^2~\text{d}\Omega +
      \frac{\triangle t}{2} \int_{\Omega_{t^{n+\frac{1}{2}}}} \nabla \cdotp (\textbf{u}^* - \textbf{w}_h )|{\textbf{u}^{n+1}_h}|^2~\text{d}\Omega\\
      &= \triangle t \int_{\Omega_{t^{n+\frac{1}{2}}}} {\bf f}^{n+\frac{1}{2}}\cdotp{\textbf{u}^{n+1}_h}~\text{d}\Omega
      \end{align}
      Now, with our consideration of $\bu^*$, the previous equation reduces to :
%

      \begin{align}
       &\ds  \int_{\Omega_{t^{n+1}}} |{\textbf{u}^{n+1}_h}|^2 ~\text{d}\Omega 
      + 2\mu \triangle t \int_{\Omega_{t^{n+\frac{1}{2}}}} |\nabla {\textbf{u}^{n+1}_h}|^2 ~\text{d}\Omega \nonumber \\ &+
      \mu_T \triangle t \int_{\Omega_{t^{n+\frac{1}{2}}}} |\nabla {\tilde{\textbf{u}}^{n+1}_h}|^2~\text{d}\Omega
      - \frac{\triangle t}{2} \int_{\Omega_{t^{n+\frac{1}{2}}}} \nabla \cdotp  \textbf{w}_h |{\textbf{u}^{n+1}_h}|^2~\text{d}\Omega\\\nonumber
      &= \triangle t \int_{\Omega_{t^{n+\frac{1}{2}}}} {\bf f}^{n+\frac{1}{2}}.{\textbf{u}^{n+1}_h}~\text{d}\Omega +
      \ds  \int_{\Omega_{t^{n}}} {\textbf{u}^{n}_h} \cdotp {\textbf{u}^{n+1}_h}~\text{d}\Omega.
      \end{align}
      Now using Holder's and  Young's  inequality we get , 
      \begin{align}
       &\ds  \int_{\Omega_{t^{n+1}}} |{\textbf{u}^{n+1}_h}|^2 ~\text{d}\Omega
      + 2\mu \triangle t \int_{\Omega_{t^{n+\frac{1}{2}}}} |\nabla {\textbf{u}^{n+1}_h}|^2 ~\text{d}\Omega \nonumber \\ \nonumber
      & +\mu_T \triangle t \int_{\Omega_{t^{n+\frac{1}{2}}}} |\nabla {\tilde{\textbf{u}}^{n+1}_h}|^2~\text{d}\Omega 
      - \frac{\triangle t}{2} \int_{\Omega_{t^{n+\frac{1}{2}}}} \nabla \cdotp  \textbf{w}_h |{\textbf{u}^{n+1}_h}|^2~\text{d}\Omega \\&\leqslant 
      \frac{1}{2} {{\parallel {\textbf{u}^{n+1}_h} \parallel}^2}_n + \frac{1}{2} {{\parallel {\textbf{u}^{n}_h} \parallel}^2}_n+
      \triangle t \frac{\mu}{2} {{\parallel {\textbf{u}^{n+1}_h} \parallel}^2}_{n+\frac{1}{2}} \\\nonumber &+
      \triangle t \Big( \frac{1+C_\Omega}{2 \mu} \Big) {{\parallel {\bf f}^{n+\frac{1}{2}} \parallel}^2}_{H^{-1}(\Omega_{t^{n+\frac{1}{2}}})},
      \end{align}
      Now breaking the first term in LHS as sum of two halves and bringing the first term of the RHS to LHS and further using Poincare's inequality we get, 
      \begin{align}
       &\ds  \frac{1}{2}\int_{\Omega_{t^{n+1}}} |{\textbf{u}^{n+1}_h}|^2~\text{d}\Omega
      + \frac{1}{2}\int_{\Omega_{t^{n+1}}} |{\textbf{u}^{n+1}_h}|^2~\text{d}\Omega
      - \frac{1}{2}\int_{\Omega_{t^{n}}} |{\textbf{u}^{n+1}_h}|^2 ~\text{d}\Omega \nonumber \\
      &+ \frac{3}{2}\mu \triangle t \int_{\Omega_{t^{n+\frac{1}{2}}}} |\nabla {\textbf{u}^{n+1}_h}|^2 ~\text{d}\Omega \label{putGCL} 
      + \mu_T \triangle t \int_{\Omega_{t^{n+\frac{1}{2}}}} |\nabla {\tilde{\textbf{u}}^{n+1}_h}|^2 ~\text{d}\Omega \\\nonumber
      &- \frac{\triangle t}{2} \int_{\Omega_{t^{n+\frac{1}{2}}}} \nabla \cdotp  \textbf{w}_h |{\textbf{u}^{n+1}_h}|^2 ~\text{d}\Omega \\\nonumber &\leqslant
	\frac{1}{2} {{\parallel {\textbf{u}^{n}_h} \parallel}^2}_n +
      \triangle t \Big( \frac{1+C_\Omega}{2 \mu} \Big) {{\parallel {\bf f}^{n+\frac{1}{2}} \parallel}^2}_{H^{-1}(\Omega_t^{n+\frac{1}{2}})}.
      \end{align}
    Further, using $\bm{\phi_i}=\bm{\phi_j}=\bu^{n+1}_h$ in GCL, \eqref{GCLcond1} we get,
    
    $$
    \frac{1}{2}\int_{\Omega_{t^{n+1}}} |{\textbf{u}^{n+1}_h}|^2~\text{d}\Omega -
    \frac{1}{2}\int_{\Omega_{t^{n}}} |{\textbf{u}^{n+1}_h}|^2~\text{d}\Omega = 
    \frac{\triangle t}{2} \int_{\Omega_{t^{n+\frac{1}{2}}}} \nabla \cdotp  \textbf{w}_h |{\textbf{u}^{n+1}_h}|^2~\text{d}\Omega.
    $$
    Now using this,  the equation \eqref{putGCL} can be reduced to,
    \begin{align}
    &\ds  \frac{1}{2}\int_{\Omega_{t^{n+1}}} |{\textbf{u}^{n+1}_h}|^2  ~\text{d}\Omega
    + \frac{3}{2}\mu \triangle t \int_{\Omega_{t^{n+\frac{1}{2}}}} |\nabla {\textbf{u}^{n+1}_h}|^2 ~\text{d}\Omega 
    + \mu_T \triangle t \int_{\Omega_{t^{n+\frac{1}{2}}}} |\nabla {\tilde{\textbf{u}}^{n+1}_h}|^2 ~\text{d}\Omega  \nonumber \\
    &\leqslant \frac{1}{2} {{\parallel {\textbf{u}^{n}_h} \parallel}^2}_n~\text{d}\Omega +
    \triangle t \Big( \frac{1+C_\Omega}{2 \mu} \Big) {{\parallel {\bf f}^{n+\frac{1}{2}} \parallel}^2}_{H^{-1}(\Omega_{t^{n+\frac{1}{2}}})}~\text{d}\Omega.
    \end{align}
    So now multiplying both sides by 2 we get,
    \begin{align}
    \ds  &\int_{\Omega_{t^{n+1}}} |{\textbf{u}^{n+1}_h}|^2  ~\text{d}\Omega
    + 3\mu \triangle t \int_{\Omega_{t^{n+\frac{1}{2}}}} |\nabla {\textbf{u}^{n+1}_h}|^2 ~\text{d}\Omega+
    \mu_T \triangle t \int_{\Omega_{t^{n+\frac{1}{2}}}} |\nabla {\tilde{\textbf{u}}^{n+1}_h}|^2 ~\text{d}\Omega+\nonumber \\
    &\leqslant  {{\parallel {\textbf{u}^{n}_h} \parallel}^2}_n+
    \triangle t \Big( \frac{1+C_\Omega}{\mu} \Big) {{\parallel {\bf f}^{n+\frac{1}{2}} \parallel}^2}_{H^{-1}(\Omega_{t^{n+\frac{1}{2}}})}.
    \end{align}
    finally summing over all time steps we obtain :
    \begin{align}\label{stabGCL}    
    \ds &\int_{\Omega_{t^{n+1}}} |{\textbf{u}^{n+1}_h}|^2  ~\text{d}\Omega
    +   \triangle t \sum\limits_{i=0}^n \Bigg( 3\mu \int_{\Omega_{t^{i+\frac{1}{2}}}} |\nabla {\textbf{u}^{i+1}_h}|^2 ~\text{d}\Omega+ 
    \mu_T\int_{\Omega_{t^{i+\frac{1}{2}}}}  |\nabla {\tilde{\textbf{u}}^{i+1}_h}|^2 ~\text{d}\Omega\Bigg)\\ \nonumber
    &\leqslant  {{\parallel {\textbf{u}^{0}_h} \parallel}^2}_0 +
    \triangle t \Big( \frac{1+C_\Omega}{\mu} \Big)\sum\limits_{i=0}^n  {{\parallel {\bf f}^{i+\frac{1}{2}} \parallel}^2}_{H^{-1}(\Omega_{t^{i+\frac{1}{2}}})}    
    \end{align}
    
    Hence, proving the result.
     \end{proof}

      \end{lemma}

      \noindent \textbf{Remark 1}: \emph{It should be noted here that the stability result we have obtained is independent
      of the domain movement or the mesh velocity }

                \section{Stability analysis of implicit Euler time discretization without GCL}
		In the previous section we derived a stability estimate for ALE-NSE in VMS form for implicit Euler time discretization scheme, and 
		because we used GCL, we arrived at an estimate which was unconditionally stable $i.e.$, independent of the domain deformation.
		Now in this section we shall derive an estimate without the use of GCL $i.e.$, by choosing a backward in time integration instead 
		of midpoint rule as done previously. This leads to a conditionally stable scheme, where the estimate depend on the mesh velocity,
		However, by choosing the time step length accordingly an estimate of the form \eqref{stabGCL}  can be obtained,
		which is unconditionally stable for all practical purposes.
		
		\begin{lemma}\label{lemma2}
		Let ${\textbf{u}^{n}_h} \in V_{h}(\Omega_{t^n})$ be as defined in the previous sections
		and the regularity assumptions
		$\textbf{ f} \in [H^{-1}(\Omega_{t})]^d$ and $\textbf{u}^{0}_h \in (L^2(\Omega))^d$ holds,
		then the following stability estimate holds for all $t \in [0,T]$
		 with the constraint  $\triangle t \leq \overline {\triangle t}$,
	       \begin{align}\label{gaststab2}
	        &\parallel \bu_h^{n+1} \parallel_{n+1}^2  + 2\mu \triangle t \sum_{i=0}^n \parallel \nabla \bu_h^{i+1} \parallel_{i+1}^2 
	        +  \mu_T \triangle t  \sum_{i=0}^n \parallel \nabla \tilde{{\textbf{u}}}^{i+1}_h\parallel_{i+1}^2\\ \nonumber
	        & \leq C' \bigg(\parallel \bu_h^{0} \parallel_{\Omega_0}^2 + \sum_{i=1}^{n+1}\parallel {\bf f} (t^i) \parallel_{H^{-1}(\Omega_{t^i})} \bigg),
	        \end{align}
	       where, $\overline {\triangle t}$ is chosen such that the following holds:
	       \begin{align}\label{boffigastcond}
	        &C (\overline {\triangle t})^2 \parallel D_\xi \hat{\textbf{w}_h} \parallel_{L^{\infty}(\Omega_0)} \parallel D_\xi \mathcal{A}_{h,  t} \parallel_{L^{\infty}(\Omega_0)} 
	       \parallel (\nabla \cdotp \textbf{w}_h ) \parallel_{L^{\infty}(\Omega_{t^{n+1}})}\\ \nonumber
	       &- \frac{\overline {\triangle t}}{2} \parallel \nabla \cdotp \bu^* \parallel_{L^{\infty}(\Omega_{t^{n+1}})} \leq \frac{1}{2}.
	       \end{align}
		
                \begin{proof}

		Let us consider the fully discrete form of ALE-NSE in VMS formulation, where for $n$= 1,2,3,...N,  find $\bu^n_h \in V_h(\Omega_t)$ such that the 
		system \eqref{ImpEulerALENSEVMS} holds for all $\bv^n_h \in V_h(\Omega_t)$,
		\begin{align}
		\ds  &\int_{\Omega_{t^{n+1}}} {\textbf{u}^{n+1}_h} \cdotp \textbf{v}_h~\text{d}\Omega  
		-  \int_{\Omega_{t^{n}}} {\textbf{u}^{n}_h} \cdotp \textbf{v}_h~\text{d}\Omega 
		+  2\mu \triangle t  \int_{\Omega_{t^{n+1}}} \nabla {\textbf{u}^{n+1}_h} \nabla \textbf{v}_h~\text{d}\Omega  \nonumber \\ \nonumber
		&+ \mu_T \triangle t  \int_{\Omega_{t^{n+1}}} \nabla \tilde{{\textbf{u}}}^{n+1}_h \nabla \tilde{\textbf{v}}_h~\text{d}\Omega +
		\triangle t \int_{\Omega_{t^{n+1}}} \nabla \cdotp [(\textbf{u}^* - \textbf{w}_h )\otimes {\textbf{u}^{n+1}_h} ] \cdotp \textbf{v}_h~\text{d}\Omega \\&= 
		\triangle t \int_{\Omega_{t^{n+1}}} {\bf f}^{n+1}.\textbf{v}_h~\text{d}\Omega
		\end{align}
	        Now taking $\bv^n_h = \bu^{n+1}_h$, we get
	        
	        \begin{align*}
	         & \ds  \int_{\Omega_{t^{n+1}}} |{\textbf{u}^{n+1}_h}|^2~\text{d}\Omega - \ds  \int_{\Omega_{t^{n}}} {\textbf{u}^{n}_h} \cdotp {\textbf{u}_h}^{n+1} ~\text{d}\Omega
	         + 2\mu \triangle t \int_{\Omega_{t^{n+1}}} |\nabla {\textbf{u}^{n+1}_h}|^2~\text{d}\Omega \nonumber\\\nonumber
                 &+ \mu_T \triangle t \int_{\Omega_{t^{n+1}}} |\nabla \tilde{{\textbf{u}}}^{n+1}_h|^2~\text{d}\Omega+
		  \triangle t \int_{\Omega_{t^{n+1}}} (\nabla \cdotp \textbf{u}^*){\textbf{u}^{n+1}_h} \cdotp{\textbf{u}^{n+1}_h}~\text{d}\Omega \nonumber \\ 
		  &+ \triangle t \int_{\Omega_{t^{n+1}}} [(\textbf{u}^*  \cdotp \nabla ){\textbf{u}^{n+1}_h}] \cdotp {\textbf{u}^{n+1}_h}~\text{d}\Omega
		  -\triangle t \int_{\Omega_{t^{n+1}}} (\nabla \cdotp \textbf{w}_h ){\textbf{u}^{n+1}_h} \cdotp{\textbf{u}^{n+1}_h}~\text{d}\Omega \nonumber \\ 
		  &- \triangle t \int_{\Omega_{t^{n+1}}} [( \textbf{w}_h  \cdotp \nabla ){\textbf{u}^{n+1}_h}]
		  \cdotp {\textbf{u}^{n+1}_h}~\text{d}\Omega		  
		  = \triangle t  \int_{\Omega_{t^{n+1}}} {\bf f}^{n+1}.{\textbf{u}^{n+1}_h}~\text{d}\Omega
	        \end{align*}
	        
	        
	        Next, combining the terms with $\bu^*$ and $\bw_h$ we get
                 
                  \begin{align}\label{staboldterm}
	         &\Rightarrow \ds  \int_{\Omega_{t^{n+1}}} |{\textbf{u}^{n+1}_h}|^2~\text{d}\Omega - \ds  \int_{\Omega_{t^{n}}} {\textbf{u}^{n}_h} \cdotp {\textbf{u}^{n+1}_h} ~\text{d}\Omega
	         + 2\mu \triangle t \int_{\Omega_{t^{n+1}}} |\nabla {\textbf{u}^{n+1}_h}|^2~\text{d}\Omega \nonumber \\
                 &+ \mu_T \triangle t \int_{\Omega_{t^{n+1}}} |\nabla \tilde{{\textbf{u}}}^{n+1}_h|^2~\text{d}\Omega
		  -\frac{\triangle t}{2} \int_{\Omega_{t^{n+1}}} (\nabla \cdotp \textbf{w}_h )| {\textbf{u}^{n+1}_h}|^2~\text{d}\Omega \\ \nonumber
		 & + \frac{\triangle t}{2}\int_{\Omega_{t^{n+1}}} (\nabla \cdotp \bu^*)|{\textbf{u}^{n+1}_h}|^2~\text{d}\Omega
		 = \triangle t  \int_{\Omega_{t^{n+1}}} {\bf f}^{n+1}\cdotp{\textbf{u}^{n+1}_h}~\text{d}\Omega.
	        \end{align}
                  In the L.H.S of the above equation we have 
	        
	         \begin{align*}
	         &\ds  \frac{1}{2}\int_{\Omega_{t^{n+1}}} |{\textbf{u}^{n+1}_h}|^2~\text{d}\Omega
	         +  \ds  \frac{1}{2}\int_{\Omega_{t^{n+1}}} |{\textbf{u}^{n+1}_h}|^2~\text{d}\Omega
	         - \ds   \frac{1}{2}\int_{\Omega_{t^{n}}} |{\textbf{u}^{n}_h}|^2  ~\text{d}\Omega \\
	        & + \ds  \frac{1}{2}\int_{\Omega_{t^{n}}} |{\textbf{u}^{n}_h} - {\textbf{u}^{n+1}_h}|^2 ~\text{d}\Omega 
	         - \ds   \frac{1}{2}\int_{\Omega_{t^{n}}} |{\textbf{u}^{n+1}_h}|^2  ~\text{d}\Omega \\
	         &+ 2\mu \triangle t \int_{\Omega_{t^{n+1}}} |\nabla {\textbf{u}^{n+1}_h}|^2~\text{d}\Omega 
                 + \mu_T \triangle t \int_{\Omega_{t^{n+1}}} |\nabla \tilde{{\textbf{u}}}^{n+1}_h|^2~\text{d}\Omega \\
		 & -\frac{\triangle t}{2} \int_{\Omega_{t^{n+1}}} (\nabla \cdotp \textbf{w}_h )| {\textbf{u}^{n+1}_h}|^2~\text{d}\Omega
		 + \frac{\triangle t}{2}\int_{\Omega_{t^{n+1}}} (\nabla \cdotp \bu^*)|{\textbf{u}^{n+1}_h}|^2~\text{d}\Omega.	 
	        \end{align*}
	        
	        Here, for the terms
	        \begin{align} \ds  \frac{1}{2}\int_{\Omega_{t^{n+1}}} |{\textbf{u}^{n+1}_h}|^2~\text{d}\Omega - \ds   \frac{1}{2}\int_{\Omega_{t^{n}}} |{\textbf{u}^{n+1}_h}|^2  ~\text{d}\Omega
	          -\frac{\triangle t}{2} \int_{\Omega_{t^{n+1}}} (\nabla \cdotp \textbf{w}_h )| {\textbf{u}^{n+1}_h}|^2~\text{d}\Omega \label{3termsGastaldiBoffi}\end{align}
	          we use Reynolds transport theorem thereby reducing it to,	          
	         \begin{align}
	          \frac{1}{2} \int_{t^n}^{t^{n+1}}\int_{\Omega_t} (\nabla \cdotp \textbf{w}_h) \textbf{u}_h^{n+1}~\text{d}\Omega dt
	          -\frac{\triangle t}{2} \int_{\Omega_{t^{n+1}}} (\nabla \cdotp \textbf{w}_h )| {\textbf{u}^{n+1}_h}|^2~\text{d}\Omega.
	         \end{align}
	          Noting that $\textbf{w}_h$ is constant in the interval $(t^n, t^{n+1})$ we can deduce that,
	           \begin{align}\label{3termsGastaldiBoffibound}
	          &  \frac{1}{2} \int_{t^n}^{t^{n+1}}\int_{\Omega_t} (\nabla \cdotp \textbf{w}_h) |\textbf{u}_h^{n+1}|^2~\text{d}\Omega dt
	          -\frac{\triangle t}{2} \int_{\Omega_{t^{n+1}}} (\nabla \cdotp \textbf{w}_h )| {\textbf{u}^{n+1}_h}|^2~\text{d}\Omega \nonumber\\
	          &= \frac{1}{2} \int_{t^n}^{t^{n+1}}\bigg(\int_{\Omega_t} (\nabla \cdotp \textbf{w}_h) |\textbf{u}_h^{n+1}|^2~\text{d}\Omega
	          - \int_{\Omega_{t^{n+1}}} (\nabla \cdotp \textbf{w}_h )| {\textbf{u}^{n+1}_h}|^2~\text{d}\Omega \bigg) dt\\ \nonumber
	          &= \frac{1}{2} \int_{t^n}^{t^{n+1}}\int_{\Omega_0} (\nabla \cdotp \textbf{w}_h) |\textbf{u}_h^{n+1}|^2 (J_{\mathcal{A}^t} - J_{\mathcal{A}^{{t}^{n+1}}}).
	           \end{align}
	          So the estimate of \eqref{3termsGastaldiBoffi} is reduced to an estimate of the difference of the Jacobians.
	          From the relations between the ALE mapping and the mesh velocity, using \eqref{Chap4DiscALE} and \eqref{Chap4Discmeshvelo},
	          we have the following estimate:
	          \begin{align}
	           |(J_{\mathcal{A}^t} - J_{\mathcal{A}^{{t}^{n+1}}})| \leq C {\triangle t} \parallel D_\xi \hat{\textbf{w}_h} \parallel_{L^{\infty}(\Omega_0)} \parallel D_\xi \mathcal{A}_{h,  t} \parallel_{L^{\infty}(\Omega_0)}.
	          \end{align}
                   Here, $D_\xi$ represents the spatial derivatives on the reference frame and C does not depend on $h, {\triangle t}$, and the ALE mapping, \cite{boffi2004stability}.
                   
                   using this estimate, \eqref{3termsGastaldiBoffibound} becomes,
	           \begin{align}\label{3termsGastaldiBoffibound1}
	          &  \frac{1}{2} \int_{t^n}^{t^{n+1}}\int_{\Omega_t} (\nabla \cdotp \textbf{w}_h) |\bv_h^{n+1}|^2~\text{d}\Omega dt
	          -\frac{\triangle t}{2} \int_{\Omega_{t^{n+1}}} (\nabla \cdotp \textbf{w}_h )| {\textbf{u}^{n+1}_h}|^2~\text{d}\Omega \nonumber \\
	          &\leq C {\triangle t}^2 \parallel D_\xi \hat{\textbf{w}_h} \parallel_{L^{\infty}(\Omega_0)} \parallel D_\xi \mathcal{A}_{h,  t} \parallel_{L^{\infty}(\Omega_0)} 
	          \int_{\Omega_{t^0}} (\nabla \cdotp \textbf{w}_h )| {\textbf{u}^{n+1}_h}|^2 J_{\mathcal{A}^{{t}^{n+1}}} ~\text{d}\Omega \\ \nonumber
	          &= C {\triangle t}^2 \parallel D_\xi \hat{\textbf{w}_h} \parallel_{L^{\infty}(\Omega_0)} \parallel D_\xi \mathcal{A}_{h,  t} \parallel_{L^{\infty}(\Omega_0)} 
	           \int_{\Omega_{t^{n+1}}} (\nabla \cdotp \textbf{w}_h )| {\textbf{u}^{n+1}_h}|^2~\text{d}\Omega
	           \end{align}
	          
	          Now using \eqref{3termsGastaldiBoffibound1}, the relation \eqref{staboldterm} can be written as the following
	          
	           \begin{align}\label{gaststab1}
	         & \ds  \frac{1}{2}\int_{\Omega_{t^{n+1}}} |{\textbf{u}^{n+1}_h}|^2~\text{d}\Omega
	          -\ds   \frac{1}{2}\int_{\Omega_{t^{n}}} |{\textbf{u}^{n}_h}|^2  ~\text{d}\Omega \nonumber\\
	         &+ 2\mu \triangle t \int_{\Omega_{t^{n+1}}} |\nabla {\textbf{u}^{n+1}_h}|^2~\text{d}\Omega
	         + \mu_T \triangle t \int_{\Omega_{t^{n+1}}} |\nabla \tilde{{\textbf{u}}}^{n+1}_h|^2~\text{d}\Omega \nonumber \\	        
	        &+ \frac{\triangle t}{2}\int_{\Omega_{t^{n+1}}} (\nabla \cdotp \bu^*)|{\textbf{u}^{n+1}_h}|^2~\text{d}\Omega\\ \nonumber
	        &- C \parallel D_\xi \hat{\textbf{w}_h} \parallel_{L^{\infty}(\Omega_0)} \parallel D_\xi \mathcal{A}_{h,  t} \parallel_{L^{\infty}(\Omega_0)}
	        {\triangle t}^2 \int_{\Omega_{t^{n+1}}} (\nabla \cdotp \textbf{w}_h )| {\textbf{u}^{n+1}_h}|^2~\text{d}\Omega\\ \nonumber
	         &\leq \triangle t  \int_{\Omega_{t^{n+1}}} {\bf f}^{n+1}\cdotp{\textbf{u}^{n+1}_h}~\text{d}\Omega.
	           \end{align}
	          Note, that this inequality is of the exact same form as equation(43) in \cite{boffi2004stability}.	       
	          Now in \eqref{gaststab1}  taking the sum over all $n$ and then using Gronwall's lemma
	          the following stability estimate can be obtained for all $\triangle t \leq \overline {\triangle t}$,
	       \begin{align}\label{gaststab2}
	        &\parallel \bu_h^{n+1} \parallel_{n+1}^2  + 2\mu \triangle t \sum_{i=0}^n \parallel \nabla \bu_h^{i+1} \parallel_{i+1}^2 
	        +  \mu_T \triangle t  \sum_{i=0}^n \parallel \nabla \tilde{{\textbf{u}}}^{i+1}_h\parallel_{i+1}^2 \nonumber \\
	        & \leq C' \bigg(\parallel \bu_h^{0} \parallel_{\Omega_0}^2 + \sum_{i=1}^{n+1}\parallel {\bf f} (t^i) \parallel_{H^{-1}(\Omega_{t^i})} \bigg),
	        \end{align}
	       where, $\overline {\triangle t}$ is chosen such that the following holds:
	       \begin{align}\label{boffigastcond}
	        &C (\overline {\triangle t})^2 \parallel D_\xi \hat{\textbf{w}_h} \parallel_{L^{\infty}(\Omega_0)} \parallel D_\xi \mathcal{A}_{h,  t} \parallel_{L^{\infty}(\Omega_0)} 
	       \parallel (\nabla \cdotp \textbf{w}_h ) \parallel_{L^{\infty}(\Omega_{t^{n+1}})}\\ \nonumber
	       &- \frac{\overline {\triangle t}}{2} \parallel \nabla \cdotp \bu^* \parallel_{L^{\infty}(\Omega_{t^{n+1}})} \leq \frac{1}{2}.
	       \end{align}
	           
                \end{proof}
	       \end{lemma}
	       Note that the stability estimate we obtained in \eqref{gaststab2} is only conditional due to the implicit dependence
	       on the mesh velocity, \eqref{boffigastcond}.
	       However, for  sufficiently small $ {\triangle t}$ we can get an estimate similar to \eqref{stabGCL}.
	       Moreover,  the LHS of  \eqref{boffigastcond} tends to zero, when $ {\triangle t} \rightarrow 0$.
	       Hence, theoretically such a choice of $\overline {\triangle t}$ is not difficult.

  Now before moving onto deriving the error estimates,  next we shall consider certain relations and inequalities.

\begin{lemma} \label{lemma3}
  The solution of the semidiscrete problem \eqref{semidiscreteALEVMS} satisfies  the following  estimates
 \begin{equation}
 \ds  \Bigg| ~\frac{d^2}{ds_2}~ \int_{\Omega_{h,s}} \bu_h(s) \cdotp  \boldsymbol{\psi} \text{d}\Omega \Bigg| \leq \text{C}(s)  \parallel \nabla \boldsymbol{\psi} \parallel_{L_{2}(\Omega_s)}
 \end{equation}
 Where C(\emph{s}) is a square integrable function on I, independent of h. 
 
 \begin{proof}
  We first define the following norm of $\bu_h$( we consider the two-dimensional case here, the three-dimensional case will follow likewise)
 \begin{equation}\label{lemma2norm1}
 \Arrowvert \ds   \bu_h  \Arrowvert_{ L^{1,2}(\Omega_t)}
=  \sum_{i=1}^{2}\bigg( \int_{\Omega_t} \vert  u_h^i \vert^2 ~\text{d}\Omega\bigg)^{\frac{1}{2}}.
\end{equation}
This is the $L_1$ norm of the vector function with $L_2$ norm of individual components, we denote it as $L^{1,2}$ norm.
The underlying function spaces being finite dimensional, this norm is equivalent to the usual
  $L_2$ norm of the vector function with $L_2$ norm of individual components, 
 \begin{equation}\label{lemma2norm2}
 \Arrowvert \ds   \bu_h  \Arrowvert_{ L^2(\Omega_t)}
=  \bigg( \sum_{i=1}^{2} \int_{\Omega_t} \vert  u_h^i \vert^2 ~\text{d}\Omega\bigg)^{\frac{1}{2}}.
 \end{equation}
We can use the definition to define $L^{1,2}$ norm for a tensor $\nabla  \bu_h$ as,
 \begin{equation}\label{lemma2norm3}
 \Arrowvert \ds \nabla  \bu_h  \Arrowvert_{ L^{1,2}(\Omega_t)}
=  \sum_{i,j=1}^{2} \bigg( \int_{\Omega_t} \bigg\vert \frac{\partial u_h^i}{\partial x_j} \bigg\vert^2 ~\text{d}\Omega\bigg)^{\frac{1}{2}} ,
 \end{equation}
 which will be equivalent to the ususal $L_2$ norm of the tensor given by,
 \begin{equation}\label{lemma2norm4}
  \Arrowvert \ds \nabla  \bu_h  \Arrowvert_{ L^{2}(\Omega_t)}
=  \bigg( \sum_{i,j=1}^{2} \int_{\Omega_t} \bigg\vert \frac{\partial u_h^i}{\partial x_j} \bigg\vert^2 ~\text{d}\Omega\bigg)^\frac{1}{2}.
 \end{equation}

An important inequality relating these two norms is 
\begin{equation}\label{normineq}
 \Arrowvert \ds   \bu_h  \Arrowvert_{ L^{1,2}(\Omega_t)} \leq \sqrt{2} \Arrowvert \ds   \bu_h  \Arrowvert_{ L^{2}(\Omega_t)}.
\end{equation}
 
      Now differentiating \eqref{semidiscreteALEVMS} w.r.t time variable t, we get 
	\begin{align}
      \ds &{\frac{d^2}{dt^2}} ( \int_{\Omega_t} \textbf{u}_h \cdotp \textbf{v}_h )~\text{d}\Omega 
      + 2\mu {\frac{d}{dt}} \int_{\Omega_t} \nabla \textbf{u}_h : \nabla \textbf{v}_h~\text{d}\Omega  
      +\mu_T {\frac{d}{dt}} \int_{\Omega_t} \nabla \tilde{\textbf{u}}_h : \nabla \tilde{\textbf{v}}_h~\text{d}\Omega \nonumber \\ \label{maineqnoflemma2}
      &+ \ds {\frac{d}{dt}} \int_{\Omega_t} \nabla \cdotp [(\textbf{u}^*- \textbf{w}_h )\otimes \textbf{u}_h ] \cdotp \textbf{v}_h~\text{d}\Omega 
      = {\frac{d}{dt}} \int_{\Omega_t} f\cdotp\textbf{v}_h~\text{d}\Omega
      \end{align}
      Here, the term on the R.H.S gives, 
	\begin{align}
      \ds {\frac{d}{dt}} \int_{\Omega_t} {\bf f} \cdotp \textbf{v}_h~\text{d}\Omega &= \sum_{i,j=1}^{2} \bigg( \ds \int_{\Omega_t} \frac{\partial f^i}{\partial t} \bigg\vert_Y  v_h^i ~\text{d}\Omega 
      +  \int_{\Omega_t} (f^i  v_h^i) (\nabla \cdotp \bw_h) ~\text{d}\Omega \bigg) \nonumber  
      \end{align}
       considering it componentwise we get  the following~\cite{NOB01}:
	\begin{align}\label{subterm01}
      &= \ds \int_{\Omega_t} \bigg( \frac{\partial f^i}{\partial t} \bigg\vert_Y ~\text{d}\Omega
      + \int_{\Omega_t}  (\nabla \cdotp \bw_h) f^i \bigg)   v_h^i  ~\text{d}\Omega \nonumber \\ \nonumber
      \\ 
      \ds &\leq \bigg(\bigg\Arrowvert \frac{\partial f^i}{\partial t} \bigg\vert_Y \bigg\Arrowvert_{ L^2(\Omega_t)} 
      + \Arrowvert (\nabla \cdotp \bw_h) \Arrowvert_{ L^\infty(\Omega_t)} \Arrowvert f^i  \Arrowvert_{ L^2(\Omega_t)} \bigg)\Arrowvert v_h^i \Arrowvert_{ L^2(\Omega_t)}
      \end{align}
      Hence, using the norms defined and \eqref{normineq}, we can write using  \eqref{subterm01}
      \begin{align}\label{term1} 
	\ds &\bigg\vert{\frac{d}{dt}} \int_{\Omega_t} {\bf f} \cdotp \textbf{v}_h~\text{d}\Omega\bigg\vert \nonumber \\ 
	&\leq \sum_{i,j=1}^{2} \bigg(\bigg\Arrowvert \frac{\partial f^i}{\partial t} \bigg\vert_Y \bigg\Arrowvert_{ L^2(\Omega_t)} 
      + \Arrowvert (\nabla \cdotp \bw_h) \Arrowvert_{ L^\infty(\Omega_t)} \Arrowvert f^i  \Arrowvert_{ L^2(\Omega_t)} \bigg)\Arrowvert v_h^i \Arrowvert_{ L^2(\Omega_t)}\\ \nonumber
      &= C_1(t) \Arrowvert \textbf{v}_h \Arrowvert_{ L^{1,2}(\Omega_t)} \leq C_2(t) \Arrowvert \textbf{v}_h \Arrowvert_{ L^{2}(\Omega_t)}
      \end{align}
      
      Now, similarly for the other terms in \eqref{maineqnoflemma2}, making use of the defined norms, 
      and the arguments presented in lemma 2.1.4  of ~\cite{NOB01}, we deduce the main result of lemma \eqref{lemma2}.

	  \end{proof}
	    \end{lemma}


	  \section{Error estimate of implicit Euler time discretization applied to conservative form of ALE-VMS}

	  In this section we shall derive an error estimate of the  implicit Euler time discretization  for a conservative form of ALE-NSE in a VMS framework.
	  We consider here the error due to the time discretization of the semidiscrete form of the equations. Our strategy is akin to \cite{NOB01}, where we manipulate the error equations
	  in a way so that the individual terms on it's LHS can be bounded using known estimates and  the stability estimates derived in the previous sections. 

	  \textbf{ Semi-discrete form of ALE-VMS (Conservative formulation) }

	  \begin{align}\label{semidisc_cons}
	  \ds &\frac{d}{dt}\int_{\Omega_t} \textbf{u}_h \cdotp \textbf{v}_h ~\text{d}\Omega_t 
	  +\int_{\Omega_t} \nabla \cdotp [(\textbf{u}^* - \textbf{w}_h) \otimes \textbf{u}_h] \cdotp \textbf{v}_h ~\text{d}\Omega_t
	  + 2\mu \int_{\Omega_t} \nabla  \textbf{u}_h : \nabla \textbf{v}_h  ~\text{d}\Omega_t  \nonumber \\
	  & + \mu_T \int_{\Omega_t} \nabla  \tilde{\textbf{u}}_h : \nabla \tilde{\textbf{v}}_h  ~\text{d}\Omega_t
	  = \int_{\Omega_t} f \cdotp \textbf{v}_h ~\text{d}\Omega_t
	  \end{align}
	  \textbf{  Fully-discrete form of ALE-VMS (Implicit Euler scheme applied to the conservative semi discrete form) }
	  \begin{align}
	  \ds &\frac{1}{\triangle t} \int_{\Omega_{t^{n+1}}} \textbf{u}^{n+1}_h \cdotp \textbf{v}_h ~\text{d}\Omega
	  - \frac{1}{\triangle t} \int_{\Omega_{t^{n}}} \textbf{u}^{n}_h \cdotp \textbf{v}_h ~\text{d}\Omega \label{impliciteuler_cons} \\ \nonumber
	  &+ \int_{\Omega_{t^{n+1}}} \nabla \cdotp [(\textbf{u}^* - \textbf{w}^{n+1}_h) \otimes \textbf{u}^{n+1}_h] \cdotp \textbf{v}_h ~\text{d}\Omega 
	  +  2\mu \int_{\Omega_{t^{n+1}}} \nabla  \textbf{u}^{n+1}_h  : \nabla \textbf{v}_h  ~\text{d}\Omega \\ \nonumber
	  & +  \ds \mu_T \int_{\Omega_{t^{n+1}}} \nabla  \tilde{\textbf{u}}^{n+1}_h : \nabla \tilde{\textbf{v}}_h  ~\text{d}\Omega
	  = \int_{\Omega_{t^{n+1}}} f^{n+1} \cdotp \textbf{v}_h ~\text{d}\Omega
	  \end{align}
	  
	  \begin{theorem}
	  The following error estimate holds for sufficiently small $\triangle t$
	      \begin{align}
	      &\parallel \textbf{u}^{n+1}_h - \textbf{u}_h(t^{n+1}) \parallel^2_{{\Omega_{t^{n+1}}}} + \triangle t \sum\limits_{i=1}^{n+1} \bigg( \mu K')\parallel \nabla \textbf{u}^{i}_h - \textbf{u}_h(t^{i}) \parallel^2_{{\Omega_{t^{i}}}}
	      + 2 \mu_T \parallel \nabla \tilde{\textbf{u}}^{i}_h - \tilde{\textbf{u}}_h(t^{i}) \parallel^2_{{\Omega_{t^{i}}}}\bigg) \nonumber \\ &\leq 
	      \mathbb{K}\triangle t^2 \Bigg[ 
	      \frac{2 }{9K} \max_{i= 1|n+1} \Big( \sup_{s\in (t^{i-1},t^{i})}\parallel J_{\mathcal{A}_{t^{i},s}} \parallel_{\text{L}^\infty(\Omega_{t^{i}})}\mathbb{Q} \Big)\\ \nonumber
	      &+ \Bigg\{ \frac{C{\triangle t}}{K}
	      \max_{i= 1|n+1} \bigg{\Arrowvert}\sup_{s\in (t^{i-1},t^{i})}\Big( \frac{\partial^2 \mathcal{A}}{\partial s_2}(Y,s) \Big)J_{\mathcal{A}_{t^{i},s}}\bigg{\Arrowvert}^2_{\text{L}^\infty(\Omega_{i})}\Bigg\}
	      \sum\limits_{i=1}^{n+1} \parallel \nabla  \textbf{u}^{i}_h \parallel^2_{\text{L}^2(\Omega_{i})}
	      \Bigg].
	      \end{align}
	     Where, $\mathbb{Q} = \ds \int_I (C(t))^2 \text{d}t$, where $C(t)$ is as defined in lemma.\ref{lemma3},  and
	     
	   $\sum\limits_{i=1}^{n+1} \parallel \nabla  \textbf{u}^{i}_h \parallel^2_{\text{L}^2(\Omega_{i})}$ is bounded by \eqref{gaststab2}.
	      
	   \begin{proof}

	 Subtracting equation \eqref{semidisc_cons} from \eqref{impliciteuler_cons} at time $t= t^{n+1}$ we get the following :
	  \begin{align*}
	  \ds &\frac{1}{\triangle t} \int_{\Omega_{t^{n+1}}} \textbf{u}^{n+1}_h \cdotp \textbf{v}_h ~\text{d}\Omega
	  - \frac{1}{\triangle t} \int_{\Omega_{t^{n}}} \textbf{u}^{n}_h \cdotp \textbf{v}_h ~\text{d}\Omega
	  - {\frac{d}{dt}\int_{\Omega_t} \textbf{u}_h \cdotp \textbf{v}_h ~\text{d}\Omega} \bigg|_{t^{n+1}} \nonumber \\ \nonumber
	  &+ 2\mu \int_{\Omega_{t^{n+1}}} \nabla  (\textbf{u}^{n+1}_h -  \textbf{u}_h(t^{n+1}) ) : \nabla  \textbf{v}_h  ~\text{d}\Omega 
	  + \ds \mu_T \int_{\Omega_{t^{n+1}}} \nabla  (\tilde{\textbf{u}}^{n+1}_h - \tilde{\textbf{u}}_h(t^{n+1}) ): \nabla \tilde{\textbf{v}}_h  ~\text{d}\Omega  \\
	  &+ \int_{\Omega_{t^{n+1}}} \nabla \cdotp (\textbf{u}^*\otimes \textbf{u}^{n+1}_h  ) \cdotp \textbf{v}_h ~\text{d}\Omega  \nonumber
	  - \int_{\Omega_{t^{n+1}}} \nabla \cdotp (\textbf{w}_h^{n+1}\otimes \textbf{u}^{n+1}_h  ) \cdotp \textbf{v}_h ~\text{d}\Omega   \\
	  &- \int_{\Omega_{t^{n+1}}} \nabla \cdotp (\textbf{u}^*\otimes \textbf{u}_h(t^{n+1})  ) \cdotp \textbf{v}_h ~\text{d}\Omega   \\ \nonumber
	  & + \int_{\Omega_{t^{n+1}}} \nabla \cdotp (\textbf{w}_h(t^{n+1})\otimes \textbf{u}_h(t^{n+1})  ) \cdotp \textbf{v}_h ~\text{d}\Omega = 0.
	  \end{align*}
	   on simplification this gives,
	  \begin{align*}
	  & \ds \frac{1}{\triangle t} \int_{\Omega_{t^{n+1}}} \textbf{u}^{n+1}_h \cdotp \textbf{v}_h ~\text{d}\Omega
	  - \frac{1}{\triangle t} \int_{\Omega_{t^{n}}} \textbf{u}^{n}_h \cdotp \textbf{v}_h ~\text{d}\Omega
	  - {\frac{d}{dt}\int_{\Omega_t} \textbf{u}_h \cdotp \textbf{v}_h ~\text{d}\Omega} \bigg|_{t^{n+1}}\nonumber \\ \nonumber
	  &+ 2\mu \int_{\Omega_{t^{n+1}}} \nabla  (\textbf{u}^{n+1}_h -  \textbf{u}_h(t^{n+1}) ) : \nabla  \textbf{v}_h  ~\text{d}\Omega 
	  + \ds \mu_T \int_{\Omega_{t^{n+1}}} \nabla  (\tilde{\textbf{u}}^{n+1}_h - \tilde{\textbf{u}}_h(t^{n+1}) ): \nabla \tilde{\textbf{v}}_h  ~\text{d}\Omega \nonumber\\
	  &+  \int_{\Omega_{t^{n+1}}} \textbf{v}_h \cdotp \nabla \cdotp [\textbf{u}^* \otimes (\textbf{u}^{n+1}_h - \textbf{u}_h(t^{n+1}))]~\text{d}\Omega \\\nonumber
	  &- \int_{\Omega_{t^{n+1}}} \textbf{v}_h \cdotp \nabla \cdotp [\textbf{w}^{n+1}_h \otimes \textbf{u}^{n+1}_h - \textbf{w}_h(t^{n+1}) \otimes \textbf{u}_h(t^{n+1})]~\text{d}\Omega = 0 
	  \end{align*}	  
	  \begin{align*}
	  &\Rightarrow \ds \frac{1}{\triangle t} \int_{\Omega_{t^{n+1}}} \textbf{u}^{n+1}_h \cdotp \textbf{v}_h ~\text{d}\Omega
	  - \frac{1}{\triangle t} \int_{\Omega_{t^{n}}} \textbf{u}^{n}_h \cdotp \textbf{v}_h ~\text{d}\Omega
	  - {\frac{d}{dt}\int_{\Omega_t} \textbf{u}_h \cdotp \textbf{v}_h ~\text{d}\Omega} \bigg|_{t^{n+1}}  \nonumber \\ \nonumber
	 & + 2\mu \int_{\Omega_{t^{n+1}}} \nabla  (\textbf{u}^{n+1}_h -  \textbf{u}_h(t^{n+1}) ): \nabla  \textbf{v}_h  ~\text{d}\Omega 
	  + \ds \mu_T \int_{\Omega_{t^{n+1}}} \nabla  (\tilde{\textbf{u}}^{n+1}_h - \tilde{\textbf{u}}_h(t^{n+1}) ): \nabla \tilde{\textbf{v}}_h  ~\text{d}\Omega \nonumber \\ 
	 & +  \int_{\Omega_{t^{n+1}}} \textbf{v}_h \cdotp \nabla \cdotp [\textbf{u}^* \otimes (\textbf{u}^{n+1}_h - \textbf{u}_h(t^{n+1}))]~\text{d}\Omega \nonumber \\
	  & - \int_{\Omega_{t^{n+1}}} \textbf{v}_h \cdotp \nabla \cdotp [\textbf{w}^{n+1}_h \otimes \textbf{u}^{n+1}_h - \textbf{w}_h(t^{n+1}) \otimes \textbf{u}_h(t^{n+1})  \\ \nonumber
	  &- \textbf{w}_h(t^{n+1}) \otimes \textbf{u}^{n+1}_h + \textbf{w}_h(t^{n+1}) \otimes \textbf{u}^{n+1}_h]~\text{d}\Omega = 0 
	  \end{align*}
	  
	  \begin{align}
	  &\Rightarrow \ds \frac{1}{\triangle t} \int_{\Omega_{t^{n+1}}} \textbf{u}^{n+1}_h \cdotp \textbf{v}_h ~\text{d}\Omega
	  - \frac{1}{\triangle t} \int_{\Omega_{t^{n}}} \textbf{u}^{n}_h \cdotp \textbf{v}_h ~\text{d}\Omega
	  - {\frac{d}{dt}\int_{\Omega_t} \textbf{u}_h \cdotp \textbf{v}_h ~\text{d}\Omega} \bigg|_{t^{n+1}}\nonumber \\ \nonumber
	  &+ 2\mu \int_{\Omega_{t^{n+1}}} \nabla  (\textbf{u}^{n+1}_h -  \textbf{u}_h(t^{n+1}) ) :\nabla  \textbf{v}_h  ~\text{d}\Omega  \label{full1}
	  + \ds \mu_T \int_{\Omega_{t^{n+1}}} \nabla  (\tilde{\textbf{u}}^{n+1}_h - \tilde{\textbf{u}}_h(t^{n+1}) ): \nabla \tilde{\textbf{v}}_h  ~\text{d}\Omega \nonumber  \\
	 & +  \int_{\Omega_{t^{n+1}}} \textbf{v}_h \cdotp \nabla \cdotp [\textbf{u}^* \otimes (\textbf{u}^{n+1}_h - \textbf{u}_h(t^{n+1}))]~\text{d}\Omega \nonumber  \\ \nonumber
	  &- \int_{\Omega_{t^{n+1}}} \textbf{v}_h \cdotp \nabla \cdotp [(\textbf{w}^{n+1}_h - \textbf{w}_h(t^{n+1})) \otimes \textbf{u}^{n+1}_h 
	  + \textbf{w}_h(t^{n+1}) \otimes (\textbf{u}^{n+1}_h - \textbf{u}_h(t^{n+1}))]~\text{d}\Omega \nonumber \\   &= 0 
	  \end{align}
	  
	  \noindent Here, we shall consider the test function $\textbf{v}_h = \textbf{u}^{n+1}_h - \textbf{u}_h(t^{n+1}) =  \textbf{e}^{n+1}_h$. 
	  Further, we manipulate terms of \eqref{full1} as follows:
	  
	  \noindent The 6th term of the \eqref{full1} gives:	  
	  \begin{align}
	  & \int_{\Omega_{t^{n+1}}} \textbf{v}_h \cdotp \nabla \cdotp [\textbf{u}^* \otimes (\textbf{u}^{n+1}_h - \textbf{u}_h(t^{n+1}))]~\text{d}\Omega\nonumber \\\nonumber
	  &= \int_{\Omega_{t^{n+1}}} \textbf{v}_h \cdotp [ \{ (\nabla \cdotp \textbf{u}^*)(\textbf{u}^{n+1}_h- \textbf{u}_h(t^{n+1}))\} 
	  + \{ (\textbf{u}^* \cdotp \nabla )(\textbf{u}^{n+1}_h - \textbf{u}_h(t^{n+1}))\} ] ~\text{d}\Omega \\ \nonumber
	  &= \int_{\Omega_{t^{n+1}}} \textbf{e}^{n+1}_h \cdotp [ \{ (\nabla \cdotp \textbf{u}^*)\textbf{e}^{n+1}_h\} 
	  + \{ (\textbf{u}^* \cdotp \nabla )\textbf{e}^{n+1}_h\} ] ~\text{d}\Omega \\ \nonumber
	  &= \int_{\Omega_{t^{n+1}}} \{ (\textbf{u}^* \cdotp \nabla) \textbf{e}^{n+1}_h \} \cdotp \textbf{e}^{n+1}_h ~\text{d}\Omega \\ \label{subterm1}
	  &= \frac{1}{2} \int_{\Omega_{t^{n+1}}} \textbf{u}^* \cdotp \{\nabla |\textbf{e}^{n+1}_h|^2 \}~\text{d}\Omega \nonumber \\
	  &= - \frac{1}{2} \int_{\Omega_{t^{n+1}}} (\nabla \cdotp \textbf{u}^*) |\textbf{e}^{n+1}_h|^2 ~\text{d}\Omega =0  
	  \end{align}
	      The 1st part of the 7th term of the \eqref{full1} gives:
	      \begin{align}
	      & \int_{\Omega_{t^{n+1}}} \textbf{v}_h \cdotp \nabla \cdotp [(\textbf{w}^{n+1}_h - \textbf{w}_h(t^{n+1})) \otimes \textbf{u}^{n+1}_h] ~\text{d}\Omega \nonumber  \\
	      &= \int_{\Omega_{t^{n+1}}} \textbf{v}_h \cdotp [\nabla \cdotp (\textbf{w}^{n+1}_h - \textbf{w}_h(t^{n+1}))\textbf{u}^{n+1}_h + 
	      \{ (\textbf{w}^{n+1}_h - \textbf{w}_h(t^{n+1}))\cdotp \nabla\}\textbf{u}^{n+1}_h] ~\text{d}\Omega \nonumber \\
	      &= \int_{\Omega_{t^{n+1}}} \nabla \cdotp (\textbf{w}^{n+1}_h - \textbf{w}_h(t^{n+1}))(\textbf{u}^{n+1}_h \cdotp \textbf{e}^{n+1}_h)  ~\text{d}\Omega \nonumber \\
		&\hspace{5mm}+ \int_{\Omega_{t^{n+1}}} [\{ (\textbf{w}^{n+1}_h - \textbf{w}_h(t^{n+1}))\cdotp \nabla\}\textbf{u}^{n+1}_h]\cdotp\textbf{e}^{n+1}_h~\text{d}\Omega \nonumber \\ 
	      &= - \int_{\Omega_{t^{n+1}}} (\textbf{w}^{n+1}_h - \textbf{w}_h(t^{n+1}))\cdotp \{\nabla(\textbf{u}^{n+1}_h \cdotp \textbf{e}^{n+1}_h)\}~\text{d}\Omega \nonumber \\
	      &\hspace{5mm}- \int_{\Omega_{t^{n+1}}} [\{ (\textbf{w}^{n+1}_h - \textbf{w}_h(t^{n+1}))\cdotp \nabla\}\textbf{e}^{n+1}_h]\cdotp\textbf{u}^{n+1}_h~\text{d}\Omega \label{subterm2}
	      \end{align}
	      In obtaining \eqref{subterm2} in the last equation we have used integration by parts and the relation $\text{b}(\textbf{u},\textbf{v},\textbf{w}) =-\text{b}(\textbf{u},\textbf{w},\textbf{v}) $,
	      where $\text{b}(\textbf{u},\textbf{v},\textbf{w}) = \int_\Omega (( \textbf{u} \cdotp \nabla) \textbf{v}) \cdotp \textbf{w}~\text{d}\Omega  $, \cite{temam1984navier} ,
	      and the 2nd part of the 7th term of the \eqref{full1} gives:
	      \begin{equation}\label{subterm3}
	      \begin{array}{rcll}
	      \begin{aligned}
	      & \int_{\Omega_{t^{n+1}}} \textbf{v}_h \cdotp \nabla \cdotp [\textbf{w}_h(t^{n+1}) \otimes (\textbf{u}^{n+1}_h - \textbf{u}_h(t^{n+1}))]~\text{d}\Omega \\
	      &= \int_{\Omega_{t^{n+1}}}  \textbf{e}^{n+1}_h  \cdotp \nabla \cdotp [\textbf{w}_h(t^{n+1}) \otimes \textbf{e}^{n+1}_h] ~\text{d}\Omega \\
	      &= \int_{\Omega_{t^{n+1}}}  \textbf{e}^{n+1}_h \cdotp \{ (\nabla  \cdotp \textbf{w}_h(t^{n+1})) \textbf{e}^{n+1}_h + (\textbf{w}_h(t^{n+1})\cdotp \nabla)\textbf{e}^{n+1}_h \}~\text{d}\Omega \\
	      &= \int_{\Omega_{t^{n+1}}} \nabla  \cdotp \textbf{w}_h(t^{n+1}) |\textbf{e}^{n+1}_h|^2 + \frac{1}{2}\textbf{w}_h(t^{n+1}) \cdotp \nabla |\textbf{e}^{n+1}_h|^2 ~\text{d}\Omega \\
	      &= \frac{1}{2}\int_{\Omega_{t^{n+1}}} \nabla  \cdotp \textbf{w}_h(t^{n+1}) |\textbf{e}^{n+1}_h|^2 ~\text{d}\Omega
	      \end{aligned}
	      \end{array}
	      \end{equation}
	      Now using \eqref{subterm1}, \eqref{subterm2}, and\eqref{subterm3}   in  \eqref{full1} 
	      we get the following:
	      \begin{align}
	      \ds &\frac{1}{\triangle t} \int_{\Omega_{t^{n+1}}} \textbf{u}^{n+1}_h \cdotp \textbf{e}^{n+1}_h ~\text{d}\Omega
	      - \frac{1}{\triangle t} \int_{\Omega_{t^{n}}} \textbf{u}^{n}_h \cdotp \textbf{e}^{n+1}_h ~\text{d}\Omega
	      - {\frac{d}{dt}\int_{\Omega_t} \textbf{u}_h \cdotp \textbf{e}^{n+1}_h~\text{d}\Omega} \bigg|_{t^{n+1}} \nonumber \\ 
	      &+ 2\mu \int_{\Omega_{t^{n+1}}} \nabla  \textbf{e}^{n+1}_h :  \nabla  \textbf{e}^{n+1}_h  ~\text{d}\Omega 
	      + \ds \mu_T \int_{\Omega_{t^{n+1}}} \nabla  \tilde{\textbf{e}}^{n+1}_h:\nabla \tilde{\textbf{e}}^{n+1}_h  ~\text{d}\Omega \nonumber \\
	     & + \int_{\Omega_{t^{n+1}}} (\textbf{w}^{n+1}_h - \textbf{w}_h(t^{n+1}))\cdotp \{\nabla(\textbf{u}^{n+1}_h \cdotp \textbf{e}^{n+1}_h)\}~\text{d}\Omega \nonumber \\ \label{full2}
	      &+ \int_{\Omega_{t^{n+1}}} [\{ (\textbf{w}^{n+1}_h - \textbf{w}_h(t^{n+1}))\cdotp \nabla\}\textbf{e}^{n+1}_h]\cdotp\textbf{u}^{n+1}_h~\text{d}\Omega \\ \nonumber
	      &- \frac{1}{2}\int_{\Omega_{t^{n+1}}} \nabla  \cdotp \textbf{w}_h(t^{n+1}) |\textbf{e}^{n+1}_h|^2 ~\text{d}\Omega =0.
	      \end{align}
	      Next, we  use the relation, due to  2.1.51 of \cite{NOB01},
	      $$\textbf{w}^{n+1}_h(\textbf{x}) - \textbf{w}_h(\textbf{x} , t^{n+1}) 
	      ~=~ - \ds \frac{1}{\triangle t} \left( \int^{t^{n+1}}_{t^n}(s- t^n) \frac{\partial^2 \mathcal{A}}{\partial s_2}(Y,s) ds\right) \circ \mathcal{A}^{-1}_{h,t^{n+1}}(\textbf{x}) $$
	      in \eqref{full2} to obtain :	 
	      \begin{align}
	      \ds &\frac{1}{\triangle t} \int_{\Omega_{t^{n+1}}} \textbf{u}^{n+1}_h \cdotp \textbf{e}^{n+1}_h ~\text{d}\Omega
	      - \frac{1}{\triangle t} \int_{\Omega_{t^{n}}} \textbf{u}^{n}_h \cdotp \textbf{e}^{n+1}_h ~\text{d}\Omega
	      - {\frac{d}{dt}\int_{\Omega_t} \textbf{u}_h \cdotp \textbf{e}^{n+1}_h~\text{d}\Omega} \bigg|_{t^{n+1}} \nonumber \\
	      &+ 2\mu \int_{\Omega_{t^{n+1}}} \nabla  \textbf{e}^{n+1}_h : \nabla  \textbf{e}^{n+1}_h  ~\text{d}\Omega 
	      + \ds \mu_T \int_{\Omega_{t^{n+1}}} \nabla  \tilde{\textbf{e}}^{n+1}_h: \nabla \tilde{\textbf{e}}^{n+1}_h  ~\text{d}\Omega \nonumber \\ &=
	      \frac{1}{\triangle t}\int_{\Omega_{t^{n+1}}} \Big\{ \Big( \int^{t^{n+1}}_{t^n}(s- t^n) \frac{\partial^2 \mathcal{A}}{\partial s_2}(Y,s) ds\Big) \circ \mathcal{A}^{-1}_{h,t^{n+1}}(\textbf{x})\Big\} 
	      \cdotp \nabla (\textbf{u}^{n+1}_h\cdotp\textbf{e}^{n+1}_h ) ~\text{d}\Omega \nonumber \\
	      &+ \int_{\Omega_{t^{n+1}}} \Big[\Big\{ \Big( \ds \frac{1}{\triangle t} \big( \int^{t^{n+1}}_{t^n}(s- t^n) \frac{\partial^2 \mathcal{A}}{\partial s_2}(Y,s) ds\big) \circ \mathcal{A}^{-1}_{h,t^{n+1}}(\textbf{x})\Big)\cdotp \nabla \Big\}
	      \textbf{e}^{n+1}_h\Big]\cdotp\textbf{u}^{n+1}_h~\text{d}\Omega \nonumber \\ \label{full3}
	      &+ \frac{1}{2}\int_{\Omega_{t^{n+1}}} \nabla  \cdotp \textbf{w}_h(t^{n+1}) |\textbf{e}^{n+1}_h|^2 ~\text{d}\Omega 
	      \end{align}

	      \noindent For ease of calculations we shall name the different terms in the R.H.S of \eqref{full3} as follows:\\

	      $\text{A} ~=~ \ds \frac{1}{\triangle t}\int_{\Omega_{t^{n+1}}} \Big\{ \Big( \int^{t^{n+1}}_{t^n}(s- t^n) \frac{\partial^2 \mathcal{A}}{\partial s_2}(Y,s) ds\Big) \circ \mathcal{A}^{-1}_{h,t^{n+1}}(\textbf{x})\Big\} 
	      \cdotp \nabla (\textbf{u}^{n+1}_h\cdotp\textbf{e}^{n+1}_h ) ~\text{d}\Omega  $ \\

	      $\text{B} ~=~ \ds \int_{\Omega_{t^{n+1}}} \Big[\Big\{ \Big( \ds \frac{1}{\triangle t} \big( \int^{t^{n+1}}_{t^n}(s- t^n) \frac{\partial^2 \mathcal{A}}{\partial s_2}(Y,s) ds\big) \circ \mathcal{A}^{-1}_{h,t^{n+1}}(\textbf{x})\Big)\cdotp \nabla \Big\}
	      \textbf{e}^{n+1}_h\Big]\cdotp\textbf{u}^{n+1}_h~\text{d}\Omega  $\\

	      $ \text{C} ~=~ \ds \frac{1}{2}\int_{\Omega_{t^{n+1}}} \nabla  \cdotp \textbf{w}_h(t^{n+1}) |\textbf{e}^{n+1}_h|^2 ~\text{d}\Omega $\\

	     \noindent Now, Using the relation 2.1.50 of \cite{NOB01}, we get:
	      \begin{equation}\label{subrelation1}
	      \begin{array}{rcll}
	      {\ds \frac{d}{dt}\int_{\Omega_t} \textbf{u}_h \cdotp \textbf{e}^{n+1}_h~\text{d}\Omega} \bigg|_{t^{n+1}} &= 
	      \ds \frac{1}{\triangle t} \int_{\Omega_{t^{n+1}}} \textbf{u}_h(t^{n+1}) \cdotp \textbf{e}^{n+1}_h ~\text{d}\Omega 
	      - \frac{1}{\triangle t} \int_{\Omega_{t^{n}}} \textbf{u}_h(t^{n}) \cdotp \textbf{e}^{n+1}_h ~\text{d}\Omega \\
	      &+ \ds \frac{1}{\triangle t} \int^{t^{n+1}}_{t^n} (s- t^n) \frac{d^2}{ds_2} \Big(\int_{\Omega_s} \textbf{u}_h(s) \cdotp  \textbf{e}^{n+1}_h \Big) ds
	      \end{array}
	      \end{equation}

	      \noindent further, using \eqref{subrelation1} and the terms denoted as A, B and C in \eqref{full3} we can write:

	      \begin{align}
	      \ds &\frac{1}{\triangle t} \int_{\Omega_{t^{n+1}}} \textbf{u}^{n+1}_h \cdotp \textbf{e}^{n+1}_h ~\text{d}\Omega
	      - \frac{1}{\triangle t} \int_{\Omega_{t^{n}}} \textbf{u}^{n}_h \cdotp \textbf{e}^{n+1}_h ~\text{d}\Omega 
	      - \frac{1}{\triangle t} \int_{\Omega_{t^{n+1}}} \textbf{u}_h(t^{n+1}) \cdotp \textbf{e}^{n+1}_h ~\text{d}\Omega \nonumber \\
	      &+ \frac{1}{\triangle t} \int_{\Omega_{t^{n}}} \textbf{u}_h(t^{n}) \cdotp \textbf{e}^{n+1}_h ~\text{d}\Omega 
	      - \ds \frac{1}{\triangle t} \int^{t^{n+1}}_{t^n} (s- t^n) \frac{d^2}{ds_2} \Big(\int_{\Omega_s} \textbf{u}_h(s) \cdotp  \textbf{e}^{n+1}_h \Big) ds \nonumber \\
	      &+ 2\mu \int_{\Omega_{t^{n+1}}} \nabla  \textbf{e}^{n+1}_h  : \nabla  \textbf{e}^{n+1}_h ~\text{d}\Omega 
	      + \ds \mu_T \int_{\Omega_{t^{n+1}}} \nabla  \tilde{\textbf{e}}^{n+1}_h : \nabla  \tilde{\textbf{e}}^{n+1}_h ~\text{d}\Omega= \text{A} +\text{B}+\text{C} \nonumber \\ \label{full4}
	      \end{align}	      
	      \noindent Combining the first four terms and further noting that
	      $~ \textbf{u}^{n+1}_h - \textbf{u}_h(t^{n+1}) =  \textbf{e}^{n+1}_h ~\text{and} ~ \textbf{u}^{n}_h - \textbf{u}_h(t^{n}) =  \textbf{e}^{n}_h$	      
	      \begin{align}
	      & \ds \frac{1}{\triangle t} \int_{\Omega_{t^{n+1}}} |\textbf{e}^{n+1}_h |^2 ~\text{d}\Omega 
	      - \frac{1}{\triangle t} \int_{\Omega_{t^{n}}} \textbf{e}^{n}_h \cdotp \textbf{e}^{n+1}_h ~\text{d}\Omega \nonumber\\ 
	      &- \ds \frac{1}{\triangle t} \int^{t^{n+1}}_{t^n} (s- t^n) \frac{d^2}{ds_2} \Big(\int_{\Omega_s} \textbf{u}_h(s) \cdotp  \textbf{e}_h \Big) ds 
	      + 2\mu \int_{\Omega_{t^{n+1}}} \nabla  \textbf{e}^{n+1}_h :\nabla  \textbf{e}^{n+1}_h ~\text{d}\Omega \\ \nonumber
	      &+ \ds \mu_T \int_{\Omega_{t^{n+1}}} \nabla  \tilde{\textbf{e}}^{n+1}_h : \nabla  \tilde{\textbf{e}}^{n+1}_h~\text{d}\Omega= \text{A} +\text{B}+\text{C}  
	      \end{align}
	      Now we shall use the following relation derived by Holder's inequality and the Reynolds transport theorem:	
	      \begin{align}
	      \int_{\Omega_{t^{n}}} \textbf{e}^{n}_h \cdotp \textbf{e}^{n+1}_h ~\text{d}\Omega &\leq \frac{1}{2} \parallel \textbf{e}^{n}_h \parallel^2_{{\Omega_{t^{n}}}} 
	      + \frac{1}{2} \parallel \textbf{e}^{n+1}_h \parallel^2_{{\Omega_{t^{n}}}}\nonumber \\ \label{relation2}
	      &= \frac{1}{2} \parallel \textbf{e}^{n}_h \parallel^2_{{\Omega_{t^{n}}}}
	      + \frac{1}{2} \parallel \textbf{e}^{n+1}_h \parallel^2_{{\Omega_{t^{n+1}}}} \\ \nonumber 
	      &- \int^{t^{n+1}}_{t^n} \int_{\Omega_s} \nabla \cdotp \textbf{w}_h(s)| \textbf{e}^{n+1}_h |^2 ~\text{d}\Omega ds
	      \end{align}
	      Now using \eqref{relation2} in \eqref{full4} we get the following :
	     \begin{align}
	      &\frac{1}{2\triangle t} \parallel \textbf{e}^{n+1}_h \parallel^2_{{\Omega_{t^{n+1}}}} + 2\mu \parallel \nabla \textbf{e}^{n+1}_h \parallel^2_{{\Omega_{t^{n+1}}}}
	      +\mu_T \parallel \nabla \tilde{\textbf{e}}^{n+1}_h \parallel^2_{{\Omega_{t^{n+1}}}} \nonumber \\
	      &\leq \frac{1}{2\triangle t} \parallel \textbf{e}^{n}_h \parallel^2_{\Omega_{t^{n}}}\nonumber \\
	      &+  \frac{1}{2}\Big( \int_{\Omega_{t^{n+1}}} \nabla  \cdotp \textbf{w}_h(t^{n+1}) |\textbf{e}^{n+1}_h|^2 ~\text{d}\Omega 
	      - \frac{1}{\triangle t}  \int^{t^{n+1}}_{t^n} \int_{\Omega_s} \nabla \cdotp \textbf{w}_h(s)| \textbf{e}^{n+1}_h |^2 ~\text{d}\Omega ds \Big) \nonumber \\
	      &+ \frac{1}{\triangle t} \int^{t^{n+1}}_{t^n} (s- t^n) \frac{d^2}{ds_2} \Big(\int_{\Omega_s} \textbf{u}_h(s) \cdotp  \textbf{e}^{n+1}_h \Big) ds \nonumber \\ 
	     & + \frac{1}{\triangle t}\int_{\Omega_{t^{n+1}}} \Big\{ \Big( \int^{t^{n+1}}_{t^n}(s- t^n) \frac{\partial^2 \mathcal{A}}{\partial s_2}(Y,s) ds\Big)
	         \circ \mathcal{A}^{-1}_{h,t^{n+1}}(\textbf{x})\Big\} 
	      \cdotp \nabla (\textbf{u}^{n+1}_h\cdotp\textbf{e}^{n+1}_h ) ~\text{d}\Omega \nonumber \\ 
	     & + \int_{\Omega_{t^{n+1}}} \Big[\Big\{ \Big( \ds \frac{1}{\triangle t} \big( \int^{t^{n+1}}_{t^n}(s- t^n) \frac{\partial^2 \mathcal{A}}{\partial s_2}(Y,s) ds\big)
	      \circ \mathcal{A}^{-1}_{h,t^{n+1}}(\textbf{x})\Big)\cdotp \nabla \Big\}
	      \textbf{e}^{n+1}_h\Big]\cdotp\textbf{u}^{n+1}_h~\text{d}\Omega \label{full5}
	      \end{align}	    
	      Now we shall denote the terms in the R.H.S of \eqref{full5} as follows for treating them separately and bounding them :	    
	      \begin{align}
	      &\text{T}_2 = \frac{1}{2}\Big( \int_{\Omega_{t^{n+1}}} \nabla  \cdotp \textbf{w}_h(t^{n+1}) |\textbf{e}^{n+1}_h|^2 ~\text{d}\Omega
	      -\frac{1}{\triangle t}  \int^{t^{n+1}}_{t^n} \int_{\Omega_s} \nabla \cdotp \textbf{w}_h(s)| \textbf{e}^{n+1}_h |^2 ~\text{d}\Omega ds \Big)  \nonumber \\ \label{subrelation2}
	     &\text{T}_3 =  \frac{1}{\triangle t} \int^{t^{n+1}}_{t^n} (s- t^n) \frac{d^2}{ds_2} \Big(\int_{\Omega_s} \textbf{u}_h(s) \cdotp  \textbf{e}^{n+1}_h \Big) ds \nonumber \\
	     &\text{T}_4 =  \frac{1}{\triangle t}\int_{\Omega_{t^{n+1}}} \Big\{ \Big( \int^{t^{n+1}}_{t^n}(s- t^n)
	                     \frac{\partial^2 \mathcal{A}}{\partial s_2}(Y,s) ds\Big) \circ \mathcal{A}^{-1}_{h,t^{n+1}}(\textbf{x})\Big\} 
	      \cdotp \nabla (\textbf{u}^{n+1}_h\cdotp\textbf{e}^{n+1}_h ) ~\text{d}\Omega \nonumber \\ \nonumber 
	     &\text{T}_5 = \int_{\Omega_{t^{n+1}}} \Big[\Big\{ \Big( \ds \frac{1}{\triangle t} \big( \int^{t^{n+1}}_{t^n}(s- t^n) \frac{\partial^2 \mathcal{A}}{\partial s_2}(Y,s) ds\big)
	      \circ \mathcal{A}^{-1}_{h,t^{n+1}}(\textbf{x})\Big)\cdotp \nabla \Big\}
	      \textbf{e}^{n+1}_h\Big]\cdotp\textbf{u}^{n+1}_h~\text{d}\Omega
	      \end{align}
	      \\\\

	      \noindent Now the term $\text{T}_2$ can be treated as follows:
	      Let us define 
	      $$
	      \gamma_i = \parallel \nabla \cdotp \textbf{w}_h(t^{i})\parallel_{\text{L}^\infty(\Omega_{t^i})} 
	      + \sup_{s\in (t^{i-1},t^{i})} \parallel J_{\mathcal{A}_{t^i,s}} \nabla \cdotp \textbf{w}_h(s) \parallel_{\text{L}^\infty(\Omega_{t^i})}
	      $$
	      where $\mathcal{A}_{t^i,s} =  \mathcal{A}_{h,s} \circ \mathcal{A}_{h,t^i} $. Hence, using this we get :
	      \begin{equation}\label{T_2}
	      \text{T}_2 \leq \frac{1}{2} \gamma_{n+1} \parallel \textbf{e}^{n+1}_h \parallel^2_{\text{L}^2(\Omega_{t^{n+1}})}.
	      \end{equation}
	      
	      \noindent Next, let us evaluate the following integral which shall be used shortly:
	      \begin{align}
	      \int^{t^{n+1}}_{t^n} (s- t^n)^2 ds &= \frac{(s- t^n)^3}{3}\bigg|^{t^{n+1}}_{t^n} \label{subrelation3} = \frac{(\triangle t)^3}{3}
	      \end{align} 
	      here, we consider $\triangle t$ to be the time step length.
	      Now using \eqref{subrelation3}, Cauchy Schwarz, Young's inequality and lemma \eqref{lemma2} we get,
	      \begin{equation}\label{T_3}
	      \begin{array}{rcll}
	      \begin{aligned}
	      \text{T}_3 &\leq \frac{1}{\triangle t} \int^{t^{n+1}}_{t^n} (s- t^n) C_4(s) \parallel \nabla  \textbf{e}^{n+1}_h \parallel_{\text{L}^2(\Omega_{s})} ds \\
	      &\leq \frac{1}{\triangle t} \int^{t^{n+1}}_{t^n} \parallel J_{\mathcal{A}_{t^{n+1},s}} \parallel^{\frac{1}{2}}_{\text{L}^\infty(\Omega_{t^{n+1}})}(s- t^n) C_4(s)
	      \parallel \nabla  \textbf{e}^{n+1}_h \parallel_{\text{L}^2(\Omega_{t^{n+1}})} ds \\
	      &\leq \frac{1}{\triangle t} \Big( \int^{t^{n+1}}_{t^n}  \parallel J_{\mathcal{A}_{t^{n+1},s}} \parallel_{\text{L}^\infty(\Omega_{t^{n+1}})}C^2_4(s) \Big)^{\frac{1}{2}}\\
	      &\hspace{12mm}\Big( \int^{t^{n+1}}_{t^n}(s- t^n)^2 \parallel \nabla  \textbf{e}^{n+1}_h \parallel^2_{\text{L}^2(\Omega_{t^{n+1}})} ds \Big)^{\frac{1}{2}}\\
	      &\leq \frac{\sqrt{\triangle t}}{3} \sup_{s\in (t^{n},t^{n+1})}\parallel J_{\mathcal{A}_{t^{n+1},s}} \parallel^{\frac{1}{2}}_{\text{L}^\infty(\Omega_{t^{n+1}})}
	      \Big( \int^{t^{n+1}}_{t^n}C^2_4(s) ds \Big)^{\frac{1}{2}}\parallel \nabla  \textbf{e}^{n+1}_h \parallel_{\text{L}^2(\Omega_{t^{n+1}})} \\
	      &\leq \frac{\triangle t}{9K} \sup_{s\in (t^{n},t^{n+1})}\parallel J_{\mathcal{A}_{t^{n+1},s}} \parallel_{\text{L}^\infty(\Omega_{t^{n+1}})}\Big( \int^{t^{n+1}}_{t^n}C^2_4(s) ds \Big)
	      + \frac{K}{4}\parallel \nabla  \textbf{e}^{n+1}_h \parallel^2_{\text{L}^2(\Omega_{t^{n+1}})}\\
	      \end{aligned}
	      \end{array}
	      \end{equation}
	      
	      \noindent Now,  by using Cauchy Schwarz, Holder's and Young's inequality,
	      and using the inverse inequality with $h= \sqrt{\triangle t}$, and further considering the property of the ALE mapping $\mathcal{A} \in [W^{1,\infty}(\Omega_0 \times I)]^d$
	      we can derive the following bound:
	      \begin{equation}\label{T_4}
	      \begin{array}{rcll}
	      \begin{aligned}
	      \text{T}_4 &= \frac{1}{\triangle t} \int^{t^{n+1}}_{t^n} \int_{\Omega_{0}} J_{\mathcal{A}_{t^{n+1},s}}(s- t^n) \nabla\widehat{(\textbf{u}^{n+1}_h 
	      \cdotp \textbf{e}^{n+1}_h)}\frac{\partial^2 \mathcal{A}}{\partial s_2}(Y,s) ~\text{d}\Omega ds \\
	      &\leq \frac{1}{\triangle t} \int^{t^{n+1}}_{t^n} \parallel J_{\mathcal{A}_{t^{n+1}}}  \parallel_{\text{L}^\infty(\Omega_0 )}
	      \bigg{\Arrowvert} \frac{\partial^2 \mathcal{A}}{\partial s_2}(Y,s) \bigg{\Arrowvert}_{\text{L}^\infty(\Omega_0)} \\
	      &\hspace{25mm}(s- t^n)\int_{\Omega_{0}} \nabla\widehat{(\textbf{u}^{n+1}_h \cdotp \textbf{e}^{n+1}_h)}~\text{d}\Omega ds\\
	      &\leq K_1 \sup_{s\in (t^{n},t^{n+1})} \bigg{\Arrowvert} \frac{\partial^2 \mathcal{A}}{\partial s_2}(Y,s) \bigg{\Arrowvert}_{\text{L}^\infty(\Omega_0)}\\& 
	      \hspace{25mm}\frac{1}{\triangle t}
	      \int^{t^{n+1}}_{t^n} (s- t^n) \parallel \nabla  (\textbf{u}^{n+1}_h \cdotp \textbf{e}^{n+1}_h) \parallel_{L^1(\Omega_{n+1})} ds\\
	      &\leq  K_1 \sup_{s\in (t^{n},t^{n+1})} \bigg{\Arrowvert} \frac{\partial^2 \mathcal{A}}{\partial s_2}(Y,s) \bigg{\Arrowvert}_{\text{L}^\infty(\Omega_0)}\\&
	      \hspace{25mm} \frac{1}{\triangle t}
	      \int^{t^{n+1}}_{t^n} (s- t^n) h^{-1} \parallel  (\textbf{u}^{n+1}_h \cdotp \textbf{e}^{n+1}_h) \parallel_{L^1(\Omega_{n+1})} ds\\
	      &\leq  K_1 \sup_{s\in (t^{n},t^{n+1})} \bigg{\Arrowvert} \frac{\partial^2 \mathcal{A}}{\partial s_2}(Y,s) \bigg{\Arrowvert}_{\text{L}^\infty(\Omega_0)}\\
	      &\hspace{25mm} \frac{1}{\triangle t}\int^{t^{n+1}}_{t^n} (s- t^n) h^{-1} \parallel  \textbf{u}^{n+1}_h \parallel_{L^2(\Omega_{n+1})} 
	      \parallel \textbf{e}^{n+1}_h \parallel_{L^2(\Omega_{n+1})} ds\\
	      &\leq  K_1 \sup_{s\in (t^{n},t^{n+1})} \bigg{\Arrowvert} \frac{\partial^2 \mathcal{A}}{\partial s_2}(Y,s) \bigg{\Arrowvert}_{\text{L}^\infty(\Omega_0)}\\ &
	      \hspace{25mm} \frac{1}{\triangle t}
	      h^{-1} \parallel  \textbf{u}^{n+1}_h \parallel_{L^2(\Omega_{n+1})}   \parallel \textbf{e}^{n+1}_h \parallel_{L^2(\Omega_{n+1})} \frac{{\triangle t}^2}{2}\\
	      &\leq K_1 \frac{\sqrt{\triangle t}}{2} \sup_{s\in (t^{n},t^{n+1})} \bigg{\Arrowvert} \frac{\partial^2 \mathcal{A}}{\partial s_2}(Y,s) \bigg{\Arrowvert}_{\text{L}^\infty(\Omega_0)}
	      \parallel  \textbf{u}^{n+1}_h \parallel_{L^2(\Omega_{n+1})}   \parallel \textbf{e}^{n+1}_h \parallel_{L^2(\Omega_{n+1})}\\
	      &\leq K \parallel \nabla  \textbf{e}^{n+1}_h \parallel^2_{L^2(\Omega_{n+1})} \\ 
	      &\hspace{25mm} + K_2 \Big( \sup_{s\in (t^{n},t^{n+1})} \bigg{\Arrowvert} \frac{\partial^2 \mathcal{A}}{\partial s_2}(Y,s) \bigg{\Arrowvert}_{\text{L}^\infty(\Omega_0)} \Big)^2
	      \frac{\triangle t}{4} \parallel \nabla  \textbf{u}^{n+1}_h \parallel^2_{L^2(\Omega_{n+1})}
	      \end{aligned}
	      \end{array}
	      \end{equation}

	      \noindent Using the relation $|\text{b}(\textbf{u},\textbf{v},\textbf{w})| \leq \parallel\textbf{u}\parallel_{\text{L}^{\infty}(\Omega)} \parallel  \nabla \textbf{v}\parallel_{\text{L}^{2}(\Omega)} 
	      \parallel\textbf{w}\parallel_{\text{L}^{2}(\Omega)}$
	      where $\text{b}(\textbf{u},\textbf{v},\textbf{w}) = \int_\Omega (( \textbf{u} \cdotp \nabla) \textbf{v}) \cdotp \textbf{w}~\text{d}\Omega  $
	      and Young's inequality we get the following:
	      \begin{equation}\label{T_5}
	      \begin{array}{rcll}
	      \begin{aligned}
	      \text{T}_5 &= \int_{\Omega_{t^{n+1}}} \Big[\Big\{ \Big( \ds \frac{1}{\triangle t} \big( \int^{t^{n+1}}_{t^n}(s- t^n) \frac{\partial^2 \mathcal{A}}{\partial s_2}(Y,s) ds\big)
	      \circ \mathcal{A}^{-1}_{h,t^{n+1}}(\textbf{x})\Big)\cdotp \nabla \Big\}
	      \textbf{e}^{n+1}_h\Big]\cdotp\textbf{u}^{n+1}_h~\text{d}\Omega \\
	      &= \int_{\Omega_{{0}}}J_{\mathcal{A}_{t^{n+1},s}} \Big[\Big\{ \Big( \ds \frac{1}{\triangle t} \big( \int^{t^{n+1}}_{t^n}(s- t^n) \frac{\partial^2 \mathcal{A}}{\partial s_2}(Y,s) ds\big)
	      \Big)\cdotp \nabla \Big\}
	      \hat{\textbf{e}}^{n+1}_h\Big]\cdotp\hat{\textbf{u}}^{n+1}_h~\text{d}\Omega \\
	      &\leq \int_{\Omega_{{0}}} 
	      \Big[\Big\{ \sup_{s\in (t^{n},t^{n+1})}\Big( \frac{\partial^2 \mathcal{A}}{\partial s_2}(Y,s) \Big)J_{\mathcal{A}_{t^{n+1},s}}\Big( \ds \frac{1}{\triangle t}  \int^{t^{n+1}}_{t^n}s- t^n ds
	      \Big)\cdotp \nabla \Big\}
	      \hat{\textbf{e}}^{n+1}_h\Big]\cdotp\hat{\textbf{u}}^{n+1}_h~\text{d}\Omega \\
	      &\leq \bigg{\Arrowvert} {\triangle t}\sup_{s\in (t^{n},t^{n+1})}\Big( \frac{\partial^2 \mathcal{A}}{\partial s_2}(Y,s) \Big)J_{\mathcal{A}_{t^{n+1},s}}\bigg{\Arrowvert}_{\text{L}^\infty(\Omega_0)}
	      \parallel \nabla  \hat{\textbf{e}}^{n+1}_h \parallel_{\text{L}^2(\Omega_0)}\parallel \hat{\textbf{u}}^{n+1}_h \parallel_{\text{L}^2(\Omega_0)}\\
	      &\leq  C{\triangle t}\bigg{\Arrowvert}\sup_{s\in (t^{n},t^{n+1})}\Big( \frac{\partial^2 \mathcal{A}}{\partial s_2}(Y,s) \Big)J_{\mathcal{A}_{t^{n+1},s}}\bigg{\Arrowvert}_{\text{L}^\infty(\Omega_{n+1})}\\
	      &\hspace{13.5mm}\parallel \nabla  \textbf{e}^{n+1}_h \parallel_{\text{L}^2(\Omega_{n+1})}\parallel\nabla   \textbf{u}^{n+1}_h \parallel_{\text{L}^2(\Omega_{n+1})}\\
	      &\leq   \frac{C{\triangle t}^2}{K}\bigg{\Arrowvert}\sup_{s\in (t^{n},t^{n+1})}\Big( \frac{\partial^2 \mathcal{A}}{\partial s_2}(Y,s) \Big)J_{\mathcal{A}_{t^{n+1},s}}\bigg{\Arrowvert}^2_{\text{L}^\infty(\Omega_{n+1})}
	      \parallel \nabla  \textbf{u}^{n+1}_h \parallel^2_{\text{L}^2(\Omega_{n+1})} \\
	      &\hspace{13.5mm} + \frac{K}{2} \parallel\nabla   \textbf{e}^{n+1}_h \parallel^2_{\text{L}^2(\Omega_{n+1})}
	      \end{aligned}
	      \end{array}
	      \end{equation}

	      Now with these bounds we have obtained for $\text{T}_2,\text{T}_3,\text{T}_4 ~\text{and}~\text{T}_5 $,  we can write the inequality  \eqref{full5} as the following :
	      \begin{align*}
	      &\frac{1}{2\triangle t} \parallel \textbf{e}^{n+1}_h \parallel^2_{{\Omega_{t^{n+1}}}} + 2\mu \parallel \nabla \textbf{e}^{n+1}_h \parallel^2_{{\Omega_{t^{n+1}}}}
	      +\mu_T \parallel \nabla \tilde{\textbf{e}}^{n+1}_h \parallel^2_{{\Omega_{t^{n+1}}}} \\ &\leq 
	      \frac{1}{2\triangle t} \parallel \textbf{e}^{n}_h \parallel^2_{\Omega_{t^{n}}}
	      +\frac{1}{2} \gamma_{n+1} \parallel \textbf{e}^{n+1}_h \parallel^2_{\text{L}^2(\Omega_{t^{n+1}})} \\
	      &+\frac{\triangle t}{9K} \sup_{s\in (t^{n},t^{n+1})}\parallel J_{\mathcal{A}_{t^{n+1},s}} \parallel_{\text{L}^\infty(\Omega_{t^{n+1}})}\Big( \int^{t^{n+1}}_{t^n}C^2_4(s) ds \Big)
	      + \frac{K}{4}\parallel \nabla  \textbf{e}^{n+1}_h \parallel^2_{\text{L}^2(\Omega_{t^{n+1}})}\\
	      &+  K \parallel \nabla  \textbf{e}^{n+1}_h \parallel^2_{L^2(\Omega_{n+1})} +
	      K_2 \Big( \sup_{s\in (t^{n},t^{n+1})} \bigg{\Arrowvert} \frac{\partial^2 \mathcal{A}}{\partial s_2}(Y,s) \bigg{\Arrowvert}_{\text{L}^\infty(\Omega_0)} \Big)^2
	      \frac{\triangle t}{4} \parallel \nabla  \textbf{u}^{n+1}_h \parallel^2_{L^2(\Omega_{n+1})}\\
	      &+  \frac{C{\triangle t}^2}{K}\bigg{\Arrowvert}\sup_{s\in (t^{n},t^{n+1})}\Big( \frac{\partial^2 \mathcal{A}}{\partial s_2}(Y,s) \Big)J_{\mathcal{A}_{t^{n+1},s}}\bigg{\Arrowvert}^2_{\text{L}^\infty(\Omega_{n+1})}
	      \parallel \nabla  \textbf{u}^{n+1}_h \parallel^2_{\text{L}^2(\Omega_{n+1})} \\ 
	      &+ \frac{K}{2} \parallel\nabla   \textbf{e}^{n+1}_h \parallel^2_{\text{L}^2(\Omega_{n+1})} \\ \\
	      \end{align*}	      
	      
	      \noindent Next we club together all the terms with $\parallel \nabla  \textbf{e}^{n+1}_h \parallel^2_{\text{L}^2(\Omega_{t^{n+1}})}$, which we would subsequently balance with the L.H.S
	      
	      \begin{align*}
	      &\frac{1}{2\triangle t} \parallel \textbf{e}^{n+1}_h \parallel^2_{{\Omega_{t^{n+1}}}} + 2\mu \parallel \nabla \textbf{e}^{n+1}_h \parallel^2_{{\Omega_{t^{n+1}}}}
	      +\mu_T \parallel \nabla \tilde{\textbf{e}}^{n+1}_h \parallel^2_{{\Omega_{t^{n+1}}}} \\ &\leq
	      \frac{1}{2\triangle t} \parallel \textbf{e}^{n}_h \parallel^2_{\Omega_{t^{n}}}+
	      \frac{1}{2} \gamma_{n+1} \parallel \textbf{e}^{n+1}_h \parallel^2_{\text{L}^2(\Omega_{t^{n+1}})} \\
	      & +\frac{\triangle t}{9K} \sup_{s\in (t^{n},t^{n+1})}\parallel J_{\mathcal{A}_{t^{n+1},s}} \parallel_{\text{L}^\infty(\Omega_{t^{n+1}})}\Big( \int^{t^{n+1}}_{t^n}C^2_4(s) ds \Big)\\
	      &+K_2 \Big( \sup_{s\in (t^{n},t^{n+1})} \bigg{\Arrowvert} \frac{\partial^2 \mathcal{A}}{\partial s_2}(Y,s) \bigg{\Arrowvert}_{\text{L}^\infty(\Omega_0)} \Big)^2
	      \frac{\triangle t}{4} \parallel \nabla  \textbf{u}^{n+1}_h \parallel^2_{L^2(\Omega_{n+1})}\\
	      &+ \frac{C{\triangle t}^2}{K}\bigg{\Arrowvert}\sup_{s\in (t^{n},t^{n+1})}\Big( \frac{\partial^2 \mathcal{A}}{\partial s_2}(Y,s) \Big)J_{\mathcal{A}_{t^{n+1},s}}\bigg{\Arrowvert}^2_{\text{L}^\infty(\Omega_{n+1})}
	      \parallel \nabla  \textbf{u}^{n+1}_h \parallel^2_{\text{L}^2(\Omega_{n+1})}\\
	      &+ \frac{7 K}{4} \parallel \nabla  \textbf{e}^{n+1}_h \parallel^2_{\text{L}^2(\Omega_{t^{n+1}})}.\\
	      \end{align*}	      
	      
	      \noindent On further simplification we get,	      
	      \begin{align}\label{final1}
	      & \parallel \textbf{e}^{n+1}_h \parallel^2_{{\Omega_{t^{n+1}}}} + 4\mu \triangle t \parallel \nabla \textbf{e}^{n+1}_h \parallel^2_{{\Omega_{t^{n+1}}}}
	      +2 \mu_T  \triangle t \parallel \nabla \tilde{\textbf{e}}^{n+1}_h \parallel^2_{{\Omega_{t^{n+1}}}} \nonumber \\ &\leq 
	      \parallel \textbf{e}^{n}_h \parallel^2_{\Omega_{t^{n}}}+
	      \triangle t\gamma_{n+1} \parallel \textbf{e}^{n+1}_h \parallel^2_{\text{L}^2(\Omega_{t^{n+1}})} \nonumber \\
	      &+\frac{2 \triangle t^2}{9K} \sup_{s\in (t^{n},t^{n+1})}\parallel J_{\mathcal{A}_{t^{n+1},s}} \parallel_{\text{L}^\infty(\Omega_{t^{n+1}})}\Big( \int^{t^{n+1}}_{t^n}C^2_4(s) ds \Big)\\ \nonumber
	      &+K_2 \Big( \sup_{s\in (t^{n},t^{n+1})} \bigg{\Arrowvert} \frac{\partial^2 \mathcal{A}}{\partial s_2}(Y,s) \bigg{\Arrowvert}_{\text{L}^\infty(\Omega_0)} \Big)^2
	      \frac{\triangle t^2}{4} \parallel \nabla  \textbf{u}^{n+1}_h \parallel^2_{L^2(\Omega_{n+1})}\\ \nonumber
	      &+ \frac{C{\triangle t}^3}{K}\bigg{\Arrowvert}\sup_{s\in (t^{n},t^{n+1})}\Big( \frac{\partial^2 \mathcal{A}}{\partial s_2}(Y,s) \Big)J_{\mathcal{A}_{t^{n+1},s}}\bigg{\Arrowvert}^2_{\text{L}^\infty(\Omega_{n+1})}
	      \parallel \nabla  \textbf{u}^{n+1}_h \parallel^2_{\text{L}^2(\Omega_{n+1})}\\ \nonumber
	     & + \frac{7 K}{2} \triangle t \parallel \nabla  \textbf{e}^{n+1}_h \parallel^2_{\text{L}^2(\Omega_{t^{n+1}})}.
	      \end{align}	      
	      Now taking summation over  n  in \eqref{final1} we get 
	      \begin{align*}
	      &\parallel \textbf{e}^{n+1}_h \parallel^2_{{\Omega_{t^{n+1}}}} + \mu K' \triangle t \sum\limits_{i=1}^{n+1} \parallel \nabla \textbf{e}^{i}_h \parallel^2_{{\Omega_{t^{i}}}}
	      +2 \mu_T  \triangle t \sum\limits_{i=0}^{n+1}\parallel \nabla \tilde{\textbf{e}}^{i}_h \parallel^2_{{\Omega_{t^{i}}}} \\ &\leq 
	      \triangle t \sum\limits_{i=1}^{n+1} \gamma_{i} \parallel \textbf{e}^{i}_h \parallel^2_{\text{L}^2(\Omega_{t^{i}})} 
	      + \frac{2 \triangle t^2}{9K} \max_{i= 1|n+1} \Big( \sup_{s\in (t^{i-1},t^{i})}\parallel J_{\mathcal{A}_{t^{i},s}} \parallel_{\text{L}^\infty(\Omega_{t^{i}})}\mathbb{Q} \Big)\\
	      &+K_2 \Big( \sup_{s\in I} \bigg{\Arrowvert} \frac{\partial^2 \mathcal{A}}{\partial s_2}(Y,s) \bigg{\Arrowvert}_{\text{L}^\infty(\Omega_0)} \Big)^2
	      \frac{\triangle t^2}{4}\sum\limits_{i=1}^{n+1} \parallel \nabla  \textbf{u}^{i}_h \parallel^2_{L^2(\Omega_{i})}\\
	      &+\frac{C{\triangle t}^3}{K}
	      \max_{i= 1|n+1} \bigg{\Arrowvert}\sup_{s\in (t^{i-1},t^{i})}\Big( \frac{\partial^2 \mathcal{A}}{\partial s_2}(Y,s) \Big)J_{\mathcal{A}_{t^{i},s}}\bigg{\Arrowvert}^2_{\text{L}^\infty(\Omega_{i})}
	      \sum\limits_{i=1}^{n+1} \parallel \nabla  \textbf{u}^{i}_h \parallel^2_{\text{L}^2(\Omega_{i})}.\\
	      \end{align*}	     
	      
	      \begin{align}\label{final2}
	      &\Rightarrow \parallel \textbf{e}^{n+1}_h \parallel^2_{{\Omega_{t^{n+1}}}} + \triangle t \sum\limits_{i=1}^{n+1} \bigg( \mu K')\parallel \nabla \textbf{e}^{i}_h \parallel^2_{{\Omega_{t^{i}}}}
	      + 2 \mu_T \parallel \nabla \tilde{\textbf{e}}^{i}_h \parallel^2_{{\Omega_{t^{i}}}}\bigg)\nonumber \\ & \leq 
	      \triangle t \sum\limits_{i=1}^{n+1}  \gamma_{i} \parallel \textbf{e}^{i}_h \parallel^2_{\text{L}^2(\Omega_{t^{i}})}
	   	      +\frac{2 \triangle t^2}{9K} \max_{i= 1|n+1} \Big( \sup_{s\in (t^{i-1},t^{i})}\parallel J_{\mathcal{A}_{t^{i},s}} \parallel_{\text{L}^\infty(\Omega_{t^{i}})}\mathbb{Q} \Big) \nonumber \\
	      &+{\triangle t^2} \Bigg[ K_2 \Big( \sup_{s\in I} \bigg{\Arrowvert} \frac{\partial^2 \mathcal{A}}{\partial s_2}(Y,s) \bigg{\Arrowvert}_{\text{L}^\infty(\Omega_0)} \Big)^2  \\ \nonumber
	      &+\frac{C{\triangle t}}{K}
	      \max_{i= 1|n+1} \bigg{\Arrowvert}\sup_{s\in (t^{i-1},t^{i})}\Big( \frac{\partial^2 \mathcal{A}}{\partial s_2}(Y,s) \Big)J_{\mathcal{A}_{t^{i},s}}\bigg{\Arrowvert}^2_{\text{L}^\infty(\Omega_{i})}\Bigg].
	      \sum\limits_{i=1}^{n+1} \parallel \nabla  \textbf{u}^{i}_h \parallel^2_{\text{L}^2(\Omega_{i})}
	      \end{align}	     
	      
	      \begin{align}\label{final3}
	      &\Rightarrow \parallel \textbf{e}^{n+1}_h \parallel^2_{{\Omega_{t^{n+1}}}} + \triangle t \sum\limits_{i=1}^{n+1} \bigg( \mu K')\parallel \nabla \textbf{e}^{i}_h \parallel^2_{{\Omega_{t^{i}}}}
	      + 2 \mu_T \parallel \nabla \tilde{\textbf{e}}^{i}_h \parallel^2_{{\Omega_{t^{i}}}}\bigg)\nonumber \\ 
	      & 
	       \triangle t \sum\limits_{i=1}^{n+1}  \gamma_{i} \parallel \textbf{e}^{i}_h \parallel^2_{\text{L}^2(\Omega_{t^{i}})} \nonumber  \\
	      &+ K_2 {\triangle t^2} \Big( \sup_{s\in I} \bigg{\Arrowvert} \frac{\partial^2 \mathcal{A}}{\partial s_2}(Y,s) \bigg{\Arrowvert}_{\text{L}^\infty(\Omega_0)} \Big)^2 
	      \sum\limits_{i=1}^{n+1} \parallel \nabla  \textbf{u}^{i}_h \parallel^2_{\text{L}^2(\Omega_{i})}\\ \nonumber
	     & +\frac{2 \triangle t^2}{9K} \max_{i= 1|n+1} \Big( \sup_{s\in (t^{i-1},t^{i})}\parallel J_{\mathcal{A}_{t^{i},s}} \parallel_{\text{L}^\infty(\Omega_{t^{i}})}\mathbb{Q} \Big)\\ \nonumber
	      &+{\triangle t^2} \Bigg[\frac{C{\triangle t}}{K}
	      \max_{i= 1|n+1} \bigg{\Arrowvert}\sup_{s\in (t^{i-1},t^{i})}\Big( \frac{\partial^2 \mathcal{A}}{\partial s_2}(Y,s) \Big)J_{\mathcal{A}_{t^{i},s}}\bigg{\Arrowvert}^2_{\text{L}^\infty(\Omega_{i})}\Bigg].
	      \sum\limits_{i=1}^{n+1} \parallel \nabla  \textbf{u}^{i}_h \parallel^2_{\text{L}^2(\Omega_{i})}
	      \end{align}	     

		Now using the discrete Gronwall's inequality, shown in lemma 1, we obtain
	      \begin{align}\label{final4}
	      &\parallel \textbf{e}^{n+1}_h \parallel^2_{{\Omega_{t^{n+1}}}} + \triangle t \sum\limits_{i=1}^{n+1} \bigg( \mu K')\parallel \nabla \textbf{e}^{i}_h \parallel^2_{{\Omega_{t^{i}}}}
	      + 2 \mu_T \parallel \nabla \tilde{\textbf{e}}^{i}_h \parallel^2_{{\Omega_{t^{i}}}}\bigg) \nonumber \\ &\leq 
	      \mathbb{K}\triangle t^2 \Bigg[ 
	      \frac{2 }{9K} \max_{i= 1|n+1} \Big( \sup_{s\in (t^{i-1},t^{i})}\parallel J_{\mathcal{A}_{t^{i},s}} \parallel_{\text{L}^\infty(\Omega_{t^{i}})}\mathbb{Q} \Big)\\ \nonumber
	      &+ \Bigg\{ \frac{C{\triangle t}}{K}
	      \max_{i= 1|n+1} \bigg{\Arrowvert}\sup_{s\in (t^{i-1},t^{i})}\Big( \frac{\partial^2 \mathcal{A}}{\partial s_2}(Y,s) \Big)J_{\mathcal{A}_{t^{i},s}}\bigg{\Arrowvert}^2_{\text{L}^\infty(\Omega_{i})}\Bigg\}
	      \sum\limits_{i=1}^{n+1} \parallel \nabla  \textbf{u}^{i}_h \parallel^2_{\text{L}^2(\Omega_{i})}
	      \Bigg].
	      \end{align}
	     
	     In \eqref{final4}, the term $\sum\limits_{i=1}^{n+1} \parallel \nabla  \textbf{u}^{i}_h \parallel^2_{\text{L}^2(\Omega_{i})}$ is dominated by using equation \eqref{gaststab2}, to get us to the final 
	     error estimate.

	     	   \end{proof}

	  \end{theorem}
	     \noindent \textbf{Remark 2}: \emph{For the error estimate we derived here we used the classical
	     form of the implicit Euler time discretization scheme, which is not restrictive in following the GCL conditions. 
	       }

\section{Summary}

This paper presents  stability and error estimates of a linearized Navier-Stokes equations
also known as the Oseen equations in a projection based variational setup. First we derive two stability estimates, both using the geometric
conservations laws (GCL) and another more general one  without using GCL. The use of GCL gives an unconditionally
stable estimate where the estimate is independent of the mesh velocity,  whereas if we choose not to use  GCL,
we have to make modification to our time steps to get to an  estimate which would be independent of the mesh velocity 
for all practical purposes. Further, we use the stability results to derive  the 
 first order error estimate due to time discretization using  backward Euler time  scheme.
 The resulting estimate can be see to  depend on the mesh velocity.

\bibliographystyle{plain}
\bibliography{masterlit}

\end{document}